\documentclass[12pt,twoside]{article}
\usepackage{mathrsfs}
\usepackage{amsfonts}
\usepackage{multirow}
\usepackage[dvips]{graphicx}

\usepackage{amsmath,amssymb}
\usepackage{mathrsfs,color}

\numberwithin{figure}{section} \numberwithin{equation}{section}
\makeatletter 
\newcommand*{\QEDA}{\hfill\ensuremath{\blacksquare}}
\begin{document}
\date{}
\author{Hengfei Ding$^{1}$,
\;\;Changpin
Li$^{2}$
\\
 \small \textit{1. School of Mathematics and Statistics, Tianshui
Normal University, Tianshui 741001, China}\\
\small \textit{2. Department of Mathematics, Shanghai University,
Shanghai 200444, China }
} \vspace{0.2cm}
\title{High-order fractional-compact finite difference method for Riesz spatial telegraph equation
\thanks{The work was partially supported by the National
Natural Science Foundation of China under Grant Nos. 11372170 and 11561060,
Key Program of Shanghai Municipal Education Commission
under Grant No. 12ZZ084, the grant of ``The First-class
Discipline of Universities in Shanghai", the Scientific Research Program
 for Young Teachers of Tianshui Normal University under Grant No. TSA1405, and Tianshui
Normal University Key Construction Subject Project (Big data processing in dynamic image).}}
\maketitle \vspace{0.2 cm}
 \hrulefill

\begin{abstract}

In this paper, we establish even order fractional-compact numerical differential formulas (4th-order, \ldots, 10th-order) for Riesz derivatives
 by using the symmetrical fractional centred difference
operator. Then we apply the derived 4th-order algorithm to the Riesz spatial telegraph equation.
 We carefully study the stability and convergence by the energy method, and show that convergence orders in temporal and spatial
directions are both 4th order.
Numerical experiments are displayed which support the fractional-compact
 numerical differentiation formulas for Riesz derivatives and the Riesz spatial telegraph equation.
\\\vspace{0.1 cm}
\vspace{0.1 cm}
\end{abstract}
\hrulefill

\section{Introduction}
In recent decades, fractional calculus (i.e., fractional integrals and fractional derivatives) have attracted increasing interest
mainly due to its potential applications
in various realms of science and engineering \cite{LiZeng,MK,SKM}.

 As far as we know, there are more than six kinds of fractional derivatives.
Usually, the time-fractional derivatives are given in the
Caputo or Riemann-Liouville sense. As for the space-fractional
derivative, it is commonly defined as an operator inverse to the Riesz potential
 and is referred to as the Riesz fractional derivative \cite{ABL,SKM}.
 Because the exact solutions of the most fractional
 differential equations cannot be analytically obtained, numerical methods for these equations
  are becoming more and more necessary and important.
Generally speaking, the first step to solve fractional differential equations is how to approximate the fractional derivatives.
 For Caputo derivatives, the commonly approximation scheme is $L_1$
formula with convergence order $\mathcal{O}(\tau^{2-\alpha})$ for $0<\alpha<1$ \cite{OS} (page 140) and \cite{LH}. In \cite{LCY}, Li et. al., proposed
 several kinds of high-order numerical schemes based
on the quadratic interpolation methods, then applied some of them to the fractional ordinary differential equations.
 Cao and Xu \cite{CX} further
developed a numerical scheme with convergence order $(3+\alpha)$ for $0<\alpha<1$ and order 4 for $\alpha>1$ under suitable
smooth conditions. Recently, a $(3-\alpha)$-th order
scheme was proposed by Gao et. al., and they used it to solve time fractional differential equations,
 but they didn't give the stability and convergence analysis for the corresponding difference scheme \cite{GSZ}.
  Almost at the same time, Li et al., \cite{LWD} also independently constructed a $(3-\alpha)$-th order numerical
   approximation formula, meanwhile they constructed two unconditionally stable difference schemes for time
fractional advection-dispersion equation with convergence orders $\mathcal{O}(\tau^{3-\alpha}+h^2)$ and $\mathcal{O}(\tau^{3-\alpha}+h^4)$,
where the rigorous theoretical analysis was displayed. More detailed introduction in this respect can be found in \cite{LiZeng}.

For Riemann-Liouville derivatives, they are often discretized by the standard Gr\"{u}nwald-Letnikov formula.
Such schemes have first order accuracy if the homogeneous initial conditions are satisfied.
 However, they are not numerically stable. To overcome this defect, Meerscharet and Tadjeran firstly proposed
  the so-called shifted Gr\"{u}nwald-Letnikov formula still with first order accuracy \cite{MT}.
  In \cite{TZD}, Tian, et.al., got a class of second order difference schemes by combining different weights and shifted
   Gr\"{u}nwald-Letnikov difference operators, the corresponding
third order quasi-compact difference schemes were proposed in \cite{ZTD}. Wang and Vong also developed a second order
difference scheme and applied it to the modified fractional anomalous subdiffusion equation \cite{WV}. In addition,
fast numerical algorithms for spatial fractional partial differential equations
with the Riemann-Liouville derivatives (or the Riesz derivatives) attract researchers' interest \cite{WB,WD}.

In general, L\'{e}vy flights and L\'{e}vy walks are used to model anomalous diffusion, which is governed
by rare but extremely large jumps of diffusing particles \cite{BHK}.
In L\'{e}vy walks, the time of travel is proportional to the total path length
if the walkers travel with constant velocity. This leads to the function with respect to time for the mean-square displacement
exists and its grows much faster than the linear case. Consequently,
This property makes L\'{e}vy walks applicable for modeling superdiffusion. However, the time
evolution of L\'{e}vy flights is simpler than that of L\'{e}vy walks. Hence,
the former is better than the latter commonly used in real life applications. In the continuous limit,
L\'{e}vy flight process is described by the superdiffusion equation, which includes the first
order derivative with respect to time and fractional Riesz operator with respect to spatial
coordinates \cite{MK},
that is to say, the Riesz space-fractional derivative has been shown to be relevant
 for anomalous diffusion models. In Riesz spatial fractional differential equations,
 the order of the Riesz derivative often belongs to $(1,2)$. Such a
Riesz derivative defined in spatial domain, which well characterize the heterogeneity and
long-range actions in space, is similar to the Laplace operator $\Delta$ in the typical partial
differential equations.
For the moment, as to the numerical algorithms for Riesz derivatives, the above numerical methods can be
 applied directly to them due to
the Riesz derivatives being the special linear
 combinations of the left
 and right Riemann-Liouville derivatives. Later,
Ortigueita defined a fractional centred
difference operator deal with the Riesz derivatives \cite{O}. \c{C}elik and Duman \cite{CD} showed that the fractional centred
difference operator has second order accuracy, and they used it to solve the Riesz spatial diffusion equation
on a finite domain. Shen et al., \cite{SLAT} also applied the above operator to the
Riesz space fractional advection diffusion equation.
 Hence, to the best of our knowledge, it seems that
 the convergence orders of the existing numerical approximation formulas for Riesz derivatives are more lower and
most of them are no more than 4, except for the results are given by Ding et al. in \cite{DL1,DL2}, who
established a series of higher-order numerical algorithms and
 also applied them to the Riesz type fractional partial differential equations, where the stability and convergence
were studied in details.
 Note that the nonlocal properties of fractional operators, using the high-order numerical algorithms for space
 fractional differential equations leads to almost the same computational cost with low-order ones
 due to the corresponding
 derived matrixes are all full and have almost the same structures. However, the convergence order
 (or accuracy) is greatly improved. So the development of the high-order numerical approximation
 formulas for Riesz derivatives have attracted utmost interest in recent years.
 In this paper, we continue to study the numerical methods for Riesz derivatives whose order lies in $(1,2)$.
 We establish the 4th-order,
 6th-order, 8th-order, and 10th-order schemes for approximating Riesz derivatives
 (the 2nd-order case was displayed in \cite{CD}), and named them as the
 fractional-compact numerical differential formulas.
 They can be smoothly applied to Riesz spatial differential equations. For convenience, we only use the fractional-compact 4th-order
 algorithm to solve Riesz spatial telegraph equation. The other order approximations can be similarly applied, so are omitted here.

The reminder of the paper is constructed as follows. In Section 2, we establish even order
approximation formulas (from 4th-order to 10th-order) for Riesz derivatives
 and study some interesting properties of the coefficients
in the fractional centred
difference operator. A temporal compact and spatial fractional-compact difference scheme for the
Riesz spatial telegraph equation is proposed in Section 3,  where the stability and convergence of
 the developed difference scheme
is also analysed. Several numerical examples are given in Section 4
to confirm the convergent orders of the proposed numerical formulas. We conclude this
paper with some remarks in the last section.

\section{The fractional-compact numerical formulas for Riesz derivatives}

\quad Firstly, we introduce the definition of the Riesz derivative and the corresponding result.

 {\it{\bf Definition 1.} The $\alpha$-order Riesz derivative of function $u(x,t)$ with respect to $x$ in real line
 for $1<\alpha<2$ is defined as \cite{SKM}
$$
\begin{array}{lll}
\displaystyle \frac{\partial^\alpha
u(x,t)}{\partial{|x|^\alpha}}&=&\displaystyle-\kappa_\alpha\left
(\,_{RL}D_{-\infty,x}^\alpha+\,_{RL}D_{x,+\infty}^\alpha\right)u(x,t)\vspace{0.2 cm}\\
&=&\displaystyle-\frac{\kappa_\alpha}{\mathrm{\Gamma}(2-\alpha)}\frac{\partial^2}
{\partial{x^2}}\int_{-\infty}^{+\infty}\frac{u(\eta,t)}{|x-\eta|^{\alpha-1}}\mathrm{d}\eta,
\end{array}\eqno(1)
$$
where $\displaystyle \kappa_\alpha
=\displaystyle\frac{1}{2}\sec\left(\frac{\pi\alpha}{2}\right)$,
$\,_{RL}D_{-\infty,x}^\alpha$ and $\,_{RL}D_{x,+\infty}^\alpha$ are the left
and right Riemann-Liouville derivatives with following forms
\begin{displaymath}
\,_{RL}D_{-\infty,x}^\alpha u(x,t)
=\frac{1}{\mathrm{\Gamma}(2-\alpha)}\frac{\partial^2}
{\partial{x^2}}\int_{-\infty}^{x}\frac{u(\eta,t)}{(x-\eta)^{\alpha-1}}\mathrm{d}\eta,\vspace{0.3cm}\\
\end{displaymath}
and
\begin{displaymath}
\,_{RL}D_{x,+\infty}^\alpha u(x,t)
=\frac{1}{\mathrm{\Gamma}(2-\alpha)}\frac{\partial^2}{\partial{x^2}}
\int_{x}^{+\infty}\frac{u(\eta,t)}{(\eta-x)^{\alpha-1}}\mathrm{d}\eta.
\end{displaymath}}
Here, $\mathrm{\Gamma}(\cdot)$ is the Euler's Gamma function,
$$
\begin{array}{lll}
\displaystyle \mathrm{\Gamma}(\beta)=\int_{0}^{\infty}\exp(-t)t^{\beta-1}\mathrm{d}t,\;\;\beta>0.
\end{array}
$$
The following function $\psi$ associated with Gamma function \cite{K} will be used later,
$$
\begin{array}{lll}
\displaystyle \psi(\beta)=\frac{\mathrm{\Gamma}'(\beta)}{\mathrm{\Gamma}(\beta)}
=-\gamma+\int_{0}^{\infty}\frac{\exp(-t)-\exp(-\beta t)}{1-\exp(-t)}\mathrm{d}t,
\end{array}\eqno(2)
$$
where $\gamma=0.57721\cdots$
 is Euler's constant.

 {\it{\bf Lemma 1. }  Let $ \alpha > 0$, $u(x,t)\in
L_1({\mathrm{R}})$ with respect to $x$, then the Fourier
 transform of Riesz derivative is given by \cite{P},
$$
\begin{array}{lll}
\displaystyle \mathscr{F}_x\left\{\frac{\partial^\alpha
u(x,t)}{\partial{|x|^\alpha}};\omega\right\}
= -|\omega|^{\alpha}\hat{u}(\omega,t),
\end{array}
$$
 where
$\hat{u}(\omega,t)$ denotes the Fourier transform of function $u(x,t)$ with respect to $x$, i.e.,
$$\displaystyle\hat{u}(\omega,t)=\mathscr{F}_x\left\{u(x,t);\omega\right\} =\int_{\mathbb{R}}\exp\left(-\mathrm{i}\omega x\right)u(x,t)\mathrm{d}x,\;\;\mathrm{i}^2=-1.$$ }

We next introduce the fractional centered difference operator
proposed by Ortigueira for $\alpha>-1$ \cite{O},
$$\begin{array}{lll}
\displaystyle
\Delta_h^{\alpha}u(x,t)=\sum\limits_{k=-\infty}^
{\infty}g_k^{(\alpha)}u(x-kh,t).
\end{array}\eqno(3)
$$
It is shown that for $1 < \alpha< 2$, then one has
$$\begin{array}{lll}
\displaystyle
\lim_{h\rightarrow0}\left(-\frac{\Delta_h^{\alpha}u(x,t)}{h^{\alpha}}\right)=\lim_{h\rightarrow0}\left(-\frac{1}{h^{\alpha}}\sum\limits_{k=-\infty}^
{\infty}g_k^{(\alpha)}u(x-kh,t)\right)
=
\frac{\partial^\alpha
u(x,t)}{\partial{|x|^\alpha}}.
\end{array}$$
where $g_k^{(\alpha)}=\frac{(-1)^{k}\mathrm{\Gamma}(\alpha+1)}{\mathrm{\Gamma}\left
(\frac{\alpha}{2}-k+1\right)\mathrm{\Gamma}\left(\frac{\alpha}{2}+k+1\right)}$, $h$ is the spatial stepsize.

In fact, the above fractional centered difference operator can be regarded as the generalization of the
second order centered difference operator, i.e.,
$$\begin{array}{lll}
\displaystyle
\lim_{\alpha\rightarrow2}\Delta_h^{\alpha}u(x,t)={u(x-h,t)-2u(x,t)+u(x+h,t)}.
\end{array}$$

 In \cite{O}, it was pointed out that
 the generating function associated with
  coefficients $\displaystyle {g}_{_k}^{(\alpha)}$ in equ. (3) is $\left|2\sin\left(\frac{\omega
h}{2}\right)\right|^{\alpha}$, i.e.,
$$\begin{array}{rrr}
\displaystyle \left|2\sin\left(\frac{\omega
h}{2}\right)\right|^{\alpha}=\sum\limits_{k=-\infty}^{\infty}{g}_{_k}^{(\alpha)}\exp(\mathrm{i}k\omega
h).
\end{array}
$$

By simple calculations, one has following recursive formula,
$$\begin{array}{rrr}
\displaystyle {g}_{_0}^{(\alpha)}=\frac{\mathrm{\Gamma}(\alpha+1)}{\mathrm{\Gamma}^2\left(\frac{\alpha}{2}+1\right)}>0,\;\;
{g}_{_k}^{(\alpha)}=\left(1-\frac{\alpha+1}{\frac{\alpha}{2}+k}\right){g}_{_{k-1}}^{(\alpha)}<0,\;k=\pm1,\pm2,\ldots.
\end{array}
$$
Noticing $0<1-\frac{\alpha+1}{\frac{\alpha}{2}+k}<1$ for $k\geq2$ and
$1-\frac{\alpha+1}{\frac{\alpha}{2}+k}>1$ for $k\leq-1$, then for fixed $\alpha\in(1,2)$ one has
$$\begin{array}{rrr}
\displaystyle
{g}_{_k}^{(\alpha)}>{g}_{_{k-1}}^{(\alpha)}>\cdots>{g}_{_{k-s}}^{(\alpha)}>\cdots,\;k=2,3,\ldots,\;s\geq1,\;k-s\geq1,
\end{array}
$$
and
$$\begin{array}{rrr}
\displaystyle
{g}_{_k}^{(\alpha)}<{g}_{_{k-1}}^{(\alpha)}<\cdots<{g}_{_{k-s}}^{(\alpha)}<\cdots,\;k=-1,-2,\ldots,\;s\leq1,\;k-s\leq-2.
\end{array}
$$

\begin{figure}[!htbp]
\centering
 \includegraphics[width=11 cm]{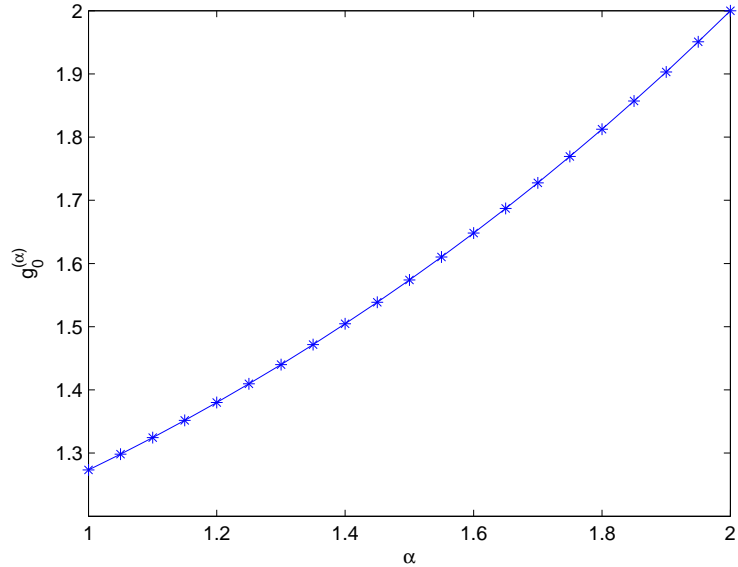}\\
  \caption{The relationship between $g_0^{(\alpha)}$ and $\alpha\in (1,2)$.}
  \label{fig.1}
\end{figure}

\begin{figure}[!htbp]
\centering
 \includegraphics[width=11 cm]{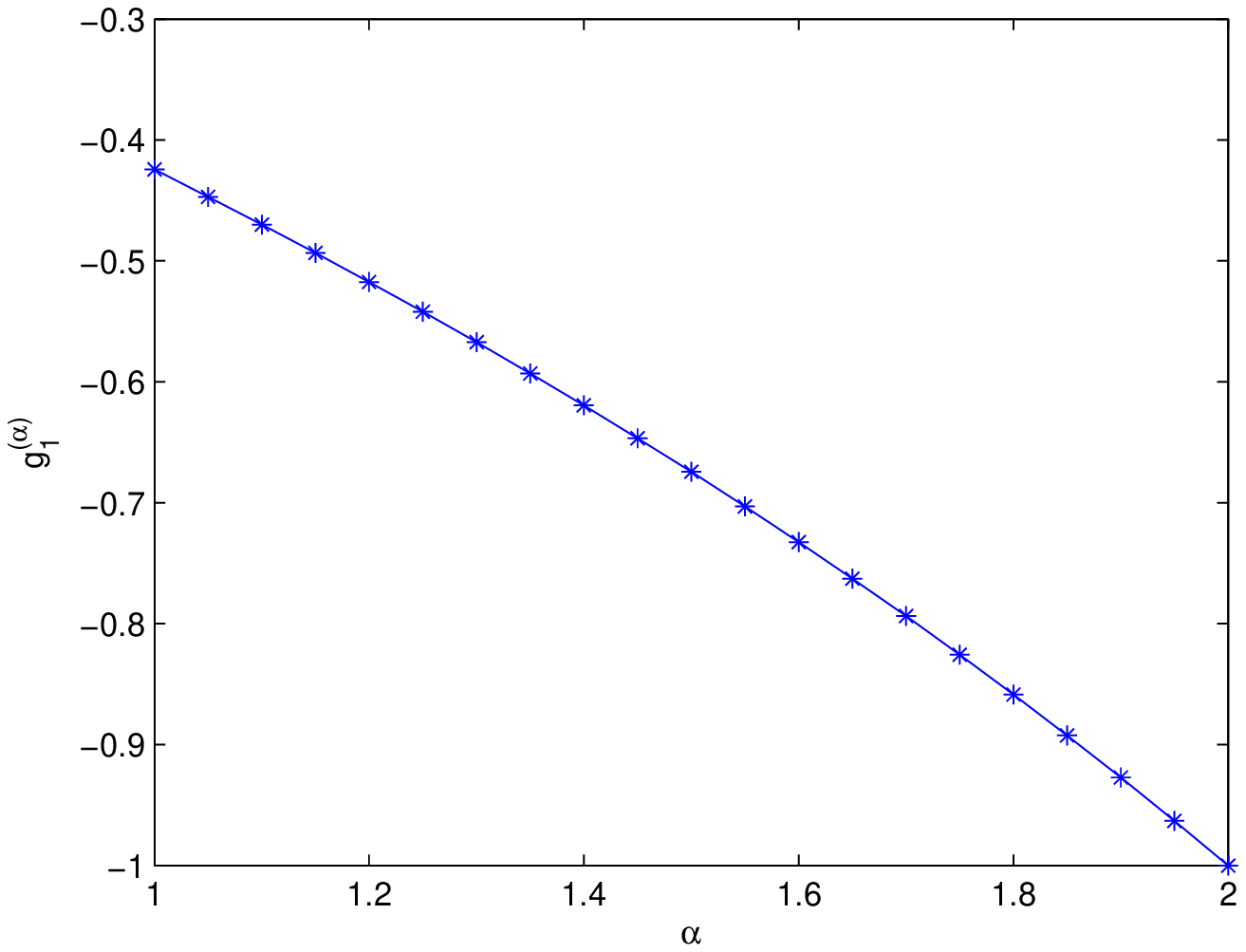}\\
  \caption{The relationship between $g_1^{(\alpha)}$ and $\alpha\in (1,2)$.}
  \label{fig.2}
\end{figure}

\begin{figure}[!htbp]
\centering
 \includegraphics[width=11 cm]{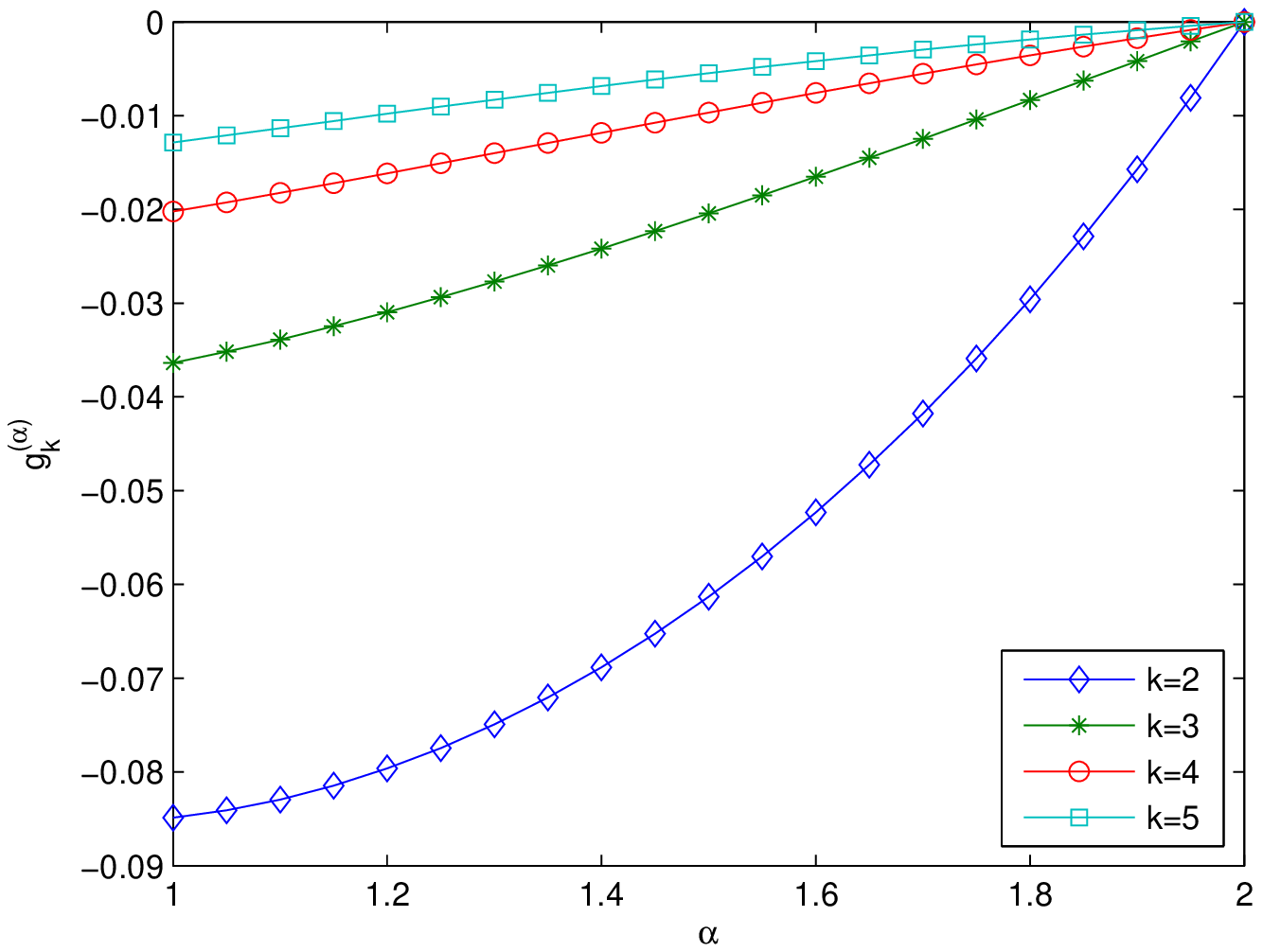}\\
  \caption{The relationship between $g_k^{(\alpha)}$ and $\alpha\in (1,2)$, $k=2,3,4,5$.}
  \label{fig.3}
\end{figure}

On the other hand, if we fix $k\geq0$, then we discuss the relationships between $\alpha$ and ${g}_{k}^{(\alpha)}$ for the
fixed $k$. Here, we plot the figures of the coefficients ${g}_{k}^{(\alpha)}$, $k=0,1,\ldots,5$.
From Figs. \ref{fig.1}--\ref{fig.3}, one can see that
${g}_{k}^{(\alpha)}$, $k=0,2,3,4,5$ increase with respect to $\alpha\in(1,2)$, ${g}_{1}^{(\alpha)}$ decrease for
$\alpha\in(1,2)$. The proof will be given later on. Next, some important properties of these coefficients are also
studied which are useful for stability and convergence analysis in some places.

{\it{\bf Theorem 1.} If $1<\alpha<2$, then the coefficients $ g_k^{(\alpha)}$ satisfy the following properties.
$$\begin{array}{lll}\displaystyle
\mathrm{(\romannumeral1)}\;\;\textit{Boundedness:}\;\; S(\alpha)\left(\frac{\alpha+4}{\alpha+2k}\right)^{2(\alpha+1)}
<\left|\right.{g}_{_k}^{(\alpha)}\left|\right.<S(\alpha)\left(\frac{\alpha+6}{\alpha+2(k+1)}\right)^{\alpha+1},\vspace{0.2 cm}\\ \;\;k\geq3, \displaystyle\;\;\textmd{where} \;\; S(\alpha)=-\frac{\mathrm{\Gamma}(\alpha+1)}
{\mathrm{\Gamma}(\frac{\alpha}{2}-1)\mathrm{\Gamma}(\frac{\alpha}{2}+3)};\vspace{0.3 cm}\\
\displaystyle
\mathrm{(\romannumeral2)}\;\; \textit{Boundedness\; of\; finite\; summations:}\vspace{0.2 cm}\\ \displaystyle
P_1(m,n,\alpha)\frac{\left(\alpha+4\right)^{2(\alpha+1)}}{2(2\alpha+1)}
<\sum\limits_{k=n}^
{m}\left|\right.{g}_{_k}^{(\alpha)}\left|\right.<P_2(m,n,\alpha)\frac{\left(\alpha+6\right)^{\alpha+1}}{2\alpha},\;\;n\geq3,\vspace{0.2 cm}\\ \displaystyle where\;\;P_1(m,n,\alpha)= -S(\alpha)\left[(\alpha+2n)^{-(2\alpha+1)}-(\alpha+2m+2)^{-(2\alpha+1)}\right],\;and\;\;\vspace{0.2 cm}\\
P_2(m,n,\alpha)= -S(\alpha)\left[(\alpha+2n)^{-\alpha}-(\alpha+2m+2)^{-\alpha}\right];
\displaystyle
\vspace{0.2 cm}\\
\displaystyle
\mathrm{(\romannumeral3)}\;\; \textit{Boundedness\; of\; infinite\; summations:}\vspace{0.2 cm}\\ \displaystyle Q_1(n,\alpha)\frac{\left(\alpha+4\right)^{2(\alpha+1)}}{2(2\alpha+1)}
<\sum\limits_{k=n}^
{\infty}\left|\right.{g}_{_k}^{(\alpha)}\left|\right.<Q_2(n,\alpha)\frac{\left(\alpha+6\right)^{\alpha+1}}{2\alpha},\;\;n\geq3,
 \vspace{0.2 cm}\\ \displaystyle where\;\;Q_1(n,\alpha)=\frac{- S(\alpha)}{(\alpha+2n)^{2\alpha+1}},\;\;and\;\;
 Q_2(n,\alpha)=\frac{- S(\alpha)}{(\alpha+2n)^{\alpha}};
\vspace{0.2 cm}\\
\displaystyle
\mathrm{(\romannumeral4)}\;\;\textit{Boundedness\; of \;two-sided\; infinite\; summations:}\;\;\\
\displaystyle
\frac{2^{1+\alpha}}{(1+\alpha)\pi}\leq\sum\limits_{k=-\infty,k\neq0}^
{\infty}\left|\right.{g}_{_k}^{(\alpha)}\left|\right.\leq\frac{2^{1+\alpha}}{\pi};
\vspace{0.2 cm}\\
\displaystyle
\mathrm{(\romannumeral5)}\;\; \textit{Monotonicity\; with\; respect\; to\;$\alpha$:}\vspace{0.2 cm}\\ \displaystyle
\textit{If}\;1\leq\alpha_1<\alpha_2\leq2, then\;{g}_{_k}^{(\alpha_1)}<{g}_{_k}^{(\alpha_2)}\;
for\; k\neq1\;and\;\;{g}_{_k}^{(\alpha_1)}>{g}_{_k}^{(\alpha_2)}\;
for\; k=1.
\end{array}
$$}
 {\bf Proof.}
See Appendix I for details.\QEDA

In the following, we construct even order fractional-compact numerical differential formulas for Riesz derivatives.

 {\it{\bf Theorem 2.} Suppose that $\displaystyle u(x,t)
\in C^{[\alpha]+2n+1}({\mathrm{R}})$ with respect to $x$ and all the derivatives
of $u(x,t)$ with respect to $x$ up to order $[\alpha]+2n+2$ exist and belong to
$\displaystyle L_{1}(\mathrm{R})$. Then the following estimate holds,}

$$\begin{array}{rrr}
\displaystyle \frac{\partial^\alpha
u(x,t)}{\partial{|x|^\alpha}}=\left(\delta_x^{0}-b_{n-1} \delta_x^{2n-2}\right)^{-1}\left(\sum\limits_{\ell=0}^
{n-2}b_{\ell}\,\delta_x^{2\ell}\right)\left(-\frac{\Delta_h^{\alpha}u(x,t)}{h^{\alpha}}\right)+\mathcal{O}(h^{2n}),\\\;n\in \mathrm{N}^{+},
\end{array}\eqno(4)
$$
{\it where
$$\begin{array}{lll}
\displaystyle
\frac{\Delta_h^{\alpha}u(x,t)}{h^{\alpha}}=\frac{1}{h^{\alpha}}\sum\limits_{k=-\infty}^
{\infty}g_k^{(\alpha)}u(x-kh,t),
\end{array}$$
and
$$\begin{array}{lll}
\displaystyle\delta_x^{2\ell}u(x_j,t)=\sum\limits_{s=0}^{2\ell}(-1)^s \dbinom{2\ell}{s} u(x_{\ell+j-s},t),\;\;\ell\geq0,
\end{array}
$$
specifically, $\delta_x^{0}$ is the identity operator with respect to spatial variable $x$, i.e., $\delta_x^{0}u(x_j,t)=u(x_j,t).$

Here, the coefficients $b_{\ell}\;\left(\ell=0,1,\ldots,n-2\right)$ satisfy the following equation
$$\begin{array}{lll}
\displaystyle
\sum\limits_{\ell=0}^
{n-2}b_{\ell}\left(2\sum\limits_{s=0}^
{\ell-1}\sum\limits_{q=0}^
{n-1}\sum\limits_{p=0}^
{n-1-q}
\frac{(-1)^{s+q}(\ell-s)^{2q}\dbinom{2\ell}{s}a_p}{(2q)!}\left|\omega h\right|^{2(p+q)}\right.\\\left.\displaystyle
+(-1)^{\ell}\dbinom{2\ell}{\ell}\sum\limits_{p=0}^
{n-1}a_p\left|\omega h\right|^{2p}\right)\\\displaystyle
=1-b_{n-1}\left(\sum\limits_{s=0}^{n-2}(-1)^s \dbinom{2n-2}{s}
\frac{2(n-1-s)^{2n-2}}{(2n-2)!}\right)(-1)^{n-1}\left|\omega h\right|^{2n-2},
\end{array}\eqno(5)
$$
and
$a_p\;(p=0,1,\ldots)$ satisfy the equation
$$\begin{array}{lll}
\displaystyle
 \sum\limits_{p=0}^
{\infty}a_p\left|\omega h\right|^{2p}&=&\displaystyle
 \left|\frac{2\sin\left(\frac{\omega h}{2}\right)}{\omega
h}\right|^{\alpha}=\left[1-\frac{\alpha}{24}|\omega
h|^2+\left(\frac{1}{1920}+\frac{\alpha-1}{1152}\right)\alpha|\omega
h|^4\right.\vspace{0.3cm}\\&&\displaystyle\left.
-\left(\frac{1}{322560}+\frac{\alpha-1}{46080}+\frac{(\alpha-1)(\alpha-2)}{82944}\right)\alpha|\omega
h|^6 +\cdots\right].
\end{array}\eqno(6)
$$
}

 {\bf Proof.}
If $n=1$, then equation (4) is reduced to the
2nd-order scheme which was reported in \cite{CD}
$$\begin{array}{lll}
\displaystyle \frac{\partial^\alpha
u(x)}{\partial{|x|^\alpha}}=-\frac{\Delta_h^{\alpha}u(x,t)}{h^\alpha}+\mathcal {O}(h^2).
\end{array}$$

The case $n=2$ was already mentioned in \cite{DL2}, here we will further
 study it and discuss its application in the next section.

 Let
$$\begin{array}{lll}
\displaystyle \mathcal {K}_n:=\left(\delta_x^{0}-b_{n-1}
\delta_x^{2n-2}\right),\;\;\; \mathcal
{L}_n:=\sum\limits_{\ell=0}^
{n-2}b_{\ell}\,\delta_x^{2\ell},\;\;n\geq2.
\end{array}$$
Taking the Fourier transform to $\mathcal {K}_n\frac{\partial^\alpha
u(x,t)}{\partial{|x|^\alpha}}$
gives
$$\begin{array}{lll}
&&\displaystyle \mathscr{F}_x\left\{\mathcal {K}_n\frac{\partial^\alpha
u(x,t)}{\partial{|x|^\alpha}};\omega\right\}=
\displaystyle
\mathscr{F}_x\left\{\left(\delta_x^{0}-b_{n-1}\delta_x^{2n-2}\right)\frac{\partial^\alpha
u(x,t)}{\partial{|x|^\alpha}};\omega\right\}
\vspace{0.2cm}\\
&=&\displaystyle-|\omega|^{\alpha}\left(1-b_{n-1}\left(\sum\limits_{s=0}^{n-2}(-1)^s \dbinom{2n-2}{s}
\frac{2(n-1-s)^{2n-2}}{(2n-2)!}(-1)^{n-1}|\omega|^{2n-2}h^{2n-2}\right.\right.\vspace{0.2 cm}\\
&&\displaystyle\left.\left.+(-1)^{n}|\omega|^{2n}\mathcal{O}(h^{2n})\right)\right)
 \widehat{u}(\omega,t),
\end{array}\eqno(7)
$$
by using Lemma 1
and the following Taylor series expansion
 $$
\begin{array}{lll}
\displaystyle
\delta_x^{2n-2}u(x,t)&=&\displaystyle \sum\limits_{s=0}^{2n-2}(-1)^s \dbinom{2n-2}{s}
u(x+(n-1-s)h,t)\vspace{0.2 cm}
\\&=&\displaystyle
\sum\limits_{s=0}^{n-2}(-1)^s \dbinom{2n-2}{s}
\left(\frac{2(n-1-s)^{2n-2}}{(2n-2)!}\frac{\partial ^{2n-2}u(x,t)}{\partial
x^{2n-2}}\right)h^{2n-2}\vspace{0.2 cm}\\&&\displaystyle
+\frac{\partial ^{2n}u(x,t)}{\partial
x^{2n}}\mathcal{O}(h^{2n}).
\end{array}
$$

Similarly, we also have
$$\begin{array}{lll}
&&\displaystyle \mathscr{F}_x\left\{\mathcal
{L}_n\left(-\frac{\Delta_h^{\alpha}u(x,t)}{h^\alpha}\right);\omega\right\}=-\displaystyle
\mathscr{F}_x\left\{\frac{1}{h^{\alpha}}\left(\sum\limits_{\ell=0}^
{n-2}b_{\ell}\,\delta_x^{2\ell}\right)\sum\limits_{k=-\infty}^
{\infty}g_k^{(\alpha)}u(x-kh);\omega\right\}
\vspace{0.2cm}\\
&=&-\displaystyle
\mathscr{F}_x\left\{\frac{1}{h^{\alpha}}\left(\sum\limits_{\ell=0}^
{n-2}b_{\ell}\sum\limits_{k=-\infty}^
{\infty}g_k^{(\alpha)}\sum\limits_{s=-\ell}^
{\ell}(-1)^{\ell-|s|}\dbinom{2\ell}{\ell-|s|}
u(x-(k+s)h)\right);\omega\right\}
\vspace{0.2cm}\\
&=&-\displaystyle\frac{1}{h^{\alpha}}\left(\sum\limits_{\ell=0}^
{n-2}b_{\ell}\sum\limits_{s=-\ell}^
{\ell}(-1)^{\ell-|s|}\dbinom{2\ell}{\ell-|s|}e^{\mathrm{i}s\omega h}\right)\sum\limits_{k=-\infty}^
{\infty}g_k^{(\alpha)}e^{\mathrm{i}k\omega h}\widehat{u}(\omega,t)
\vspace{0.2cm}\\
&=&-\displaystyle\frac{1}{h^{\alpha}}\left(\sum\limits_{\ell=0}^
{n-2}b_{\ell}\left(\sum\limits_{s=0}^
{\ell-1}2(-1)^{s}\dbinom{2\ell}{s}\cos((\ell-s)\omega h)+(-1)^{\ell}\dbinom{2\ell}{\ell}\right)\right)\vspace{0.2cm}\\&&\displaystyle\times\sum\limits_{k=-\infty}^
{\infty}g_k^{(\alpha)}e^{\mathrm{i}k\omega h}\widehat{u}(\omega,t)
\vspace{0.2cm}\\

&=&-\displaystyle\frac{1}{h^{\alpha}}\left(\sum\limits_{\ell=0}^
{n-2}b_{\ell}\left(\sum\limits_{s=0}^
{\ell-1}2(-1)^{s}\dbinom{2\ell}{s}\sum\limits_{q=0}^
{\infty}(-1)^q\frac{((\ell-s)\omega h)^{2q}}{(2q)!}+(-1)^{\ell}\dbinom{2\ell}{\ell}\right)\right)
\vspace{0.2cm}\\&&\displaystyle\times\sum\limits_{k=-\infty}^
{\infty}g_k^{(\alpha)}e^{\mathrm{i}k\omega h}\widehat{u}(\omega,t)
\vspace{0.2cm}\\

\end{array}
$$
$$\begin{array}{lll}
&=&\displaystyle-|\omega|^{\alpha}\left(\left(\sum\limits_{\ell=0}^
{n-2}b_{\ell}\left(\sum\limits_{s=0}^
{\ell-1}2(-1)^{s}\dbinom{2\ell}{s}\sum\limits_{q=0}^
{n-1}(-1)^q\frac{((\ell-s)\omega h)^{2q}}{(2q)!}\right.\right.\right.\vspace{0.2cm}\\&&\left.\left.\left.+(-1)^{\ell}\dbinom{2\ell}{\ell}\right)\right)\left|\frac{2\sin\left(\frac{\omega
h}{2}\right)}{\omega h}\right|^{\alpha}+|\omega|^{2n}\mathcal {O}(h^{2n})\right) \widehat{u}(\omega,t)\\
&=&\displaystyle-|\omega|^{\alpha}\left(\left(\sum\limits_{\ell=0}^
{n-2}b_{\ell}\left(\sum\limits_{s=0}^
{\ell-1}2(-1)^{s}\dbinom{2\ell}{s}\sum\limits_{q=0}^
{n-1}(-1)^q\frac{((\ell-s)\omega h)^{2q}}{(2q)!}\right.\right.\right.\vspace{0.2cm}\\&&\displaystyle\left.\left.\left.+(-1)^{\ell}\dbinom{2\ell}{\ell}\right)\right)\sum\limits_{p=0}^
{n-1}a_p\left|\omega h\right|^{2p}+|\omega|^{2n}\mathcal {O}(h^{2n})\right) \widehat{u}(\omega,t).
\end{array}\eqno(8)
$$
After noting equs. (7) and (8), one gets
 $$
\begin{array}{lll}
\displaystyle
\widehat\delta(\omega,h)&=&\displaystyle\mathscr{F}_x\left\{\mathcal {K}_n\frac{\partial^\alpha
u(xt,)}{\partial{|x|^\alpha}};\omega\right\}-
\mathscr{F}_x\left\{\mathcal
{L}_n\left(-\frac{\Delta_h^{\alpha}u(x,t)}{h^\alpha}\right);\omega\right\}\vspace{0.2cm}\\
&=&\displaystyle \widehat{C}_1h^{2n}|\omega|^{2n+\alpha} \widehat{u}(\omega,t).
\end{array}\eqno(9)
$$

Since $u(x,t)\in C^{[\alpha]+2n+1}(\mathbb{\mathrm{R}})$ and its partial derivatives up
to order $[\alpha]+2n+2$ with respect to $x$ belong to $L_1(\mathbb{\mathrm{R}})$,
there exists a positive constant $\widehat{C}_2$ such that
$$\begin{array}{rrr}
\displaystyle \left|\widehat{u}(\omega,t)\right|\leq
\widehat{C}_2\left(1+\left|\omega\right|\right)^{-([\alpha]+2n+2)}.
\end{array}
$$
So, using (9) we receive that
$$\begin{array}{lll}
\displaystyle \left|\widehat{\delta}(\omega,h)\right|\leq
\widehat{C}_1\left|\omega\right|^{\alpha+2n}h^{2n}\left|\widehat{u}(\omega,t)\right|&\leq&\displaystyle
\widehat{C}_1\left(1+\left|\omega\right|\right)^{\alpha+2n}h^{2n}\left|\widehat{u}(\omega,t)\right|\\&\leq&\displaystyle
\widehat{C}_3\left(1+\left|\omega\right|\right)^{\alpha-[\alpha]-2}h^{2n},
\end{array}\eqno(10)
$$
where $\widehat{C}_3=\widehat{C}_1\widehat{C}_2$.

Furthermore, taking the inverse Fourier transform on both sides of
(9) and combining with (10) give
$$
\begin{array}{lll}
\displaystyle \left|\mathcal {K}_n\frac{\partial^\alpha
u(x,t)}{\partial{|x|^\alpha}}-
\mathcal
{L}_n\left(-\frac{\Delta_h^{\alpha}u(x,t)}{h^\alpha}\right)\right|
&=&\left|\delta(\omega,h)\right|
=\displaystyle\frac{1}{2\pi}\left|\int_{\mathbb{R}}\widehat{\delta}(\omega,h)\exp(\mathrm{i}\omega
 h)\mathrm{d}\omega\right|\vspace{0.2 cm}\\ &\leq&\displaystyle
\frac{1}{2\pi}\int_{\mathbb{R}}\left|\widehat{\delta}(\omega,h)\right|\mathrm{d}\omega
\vspace{0.2 cm}\\ &\leq&\displaystyle
\frac{\widehat{C}_3}{2\pi}\left(\int_{\mathbb{R}}\left(1+\left|\omega\right|\right)^{\alpha-[\alpha]-2}\mathrm{d}\omega\right)h^{2n}
\vspace{0.2 cm}\\ &=&\displaystyle \widehat{C}h^{2n},
\end{array}$$
i.e.,
$$\begin{array}{lll}
\displaystyle \frac{\partial^\alpha
u(x,t)}{\partial{|x|^\alpha}}=\left(\delta_x^{0}-b_{n-1} \delta_x^{2n-2}\right)^{-1}\left(\sum\limits_{\ell=0}^
{n-2}b_{\ell}\,\delta_x^{2\ell}\right)\left(-\frac{\Delta_h^{\alpha}u(x,t)}{h^\alpha}\right)+\mathcal{O}(h^{2n}),\;n\in\mathrm{ N}^{+},
\end{array}
$$
where $\displaystyle
\widehat{C}=\frac{\widehat{C}_3}{(1+[\alpha]-\alpha)\pi}$. This completes the proof.\QEDA

Here, we list the common even order fractional-compact numerical differential formulas
 in view of equs. (5) and (6):
$$\begin{array}{lll}
\displaystyle \frac{\partial^\alpha
u(x,t)}{\partial{|x|^\alpha}}=\left({\delta_x^{0}-b_1\delta_x^2}\right)^{-1}b_0\delta_x^{0}\left(-\frac{\Delta_h^{\alpha}u(x,t)}{h^\alpha}\right)
+ \mathcal{O}(h^4),
\end{array}\eqno(11)
$$
$$\begin{array}{lll}
\displaystyle \frac{\partial^\alpha
u(x,t)}{\partial{|x|^\alpha}}=\left(\delta_x^{0}-b_2\delta_x^4\right)^{-1}(b_0\delta_x^{0}+b_1\delta_x^2)\left(-\frac{\Delta_h^{\alpha}u(x,t)}{h^\alpha}\right)
+ \mathcal{O}(h^6),
\end{array}\eqno(12)
$$
$$\begin{array}{lll}
\displaystyle \frac{\partial^\alpha
u(x,t)}{\partial{|x|^\alpha}}&=&\displaystyle\left(\delta_x^{0}-b_3\delta_x^6\right)^{-1}(b_0\delta_x^{0}+b_1\delta_x^2+b_2\delta_x^4)
\left(-\frac{\Delta_h^{\alpha}u(x,t)}{h^\alpha}\right)\vspace{0.2 cm}\\&&\displaystyle
+ \mathcal{O}(h^8),
\end{array}\eqno(13)
$$
$$\begin{array}{lll}
\displaystyle \frac{\partial^\alpha
u(x,t)}{\partial{|x|^\alpha}}&=&\displaystyle\left(\delta_x^{0}-b_4\delta_x^8\right)^{-1}(b_0\delta_x^{0}+b_1\delta_x^2+b_2\delta_x^4+b_3\delta_x^6)
\left(-\frac{\Delta_h^{\alpha}u(x,t)}{h^\alpha}\right) \vspace{0.2 cm}\\&&\displaystyle+ \mathcal{O}(h^{10}),
\end{array}\eqno(14)
$$
where
$$\begin{array}{lll}
\displaystyle b_0=1,\;\;\vspace{0.3 cm}\\
\displaystyle b_1=-\frac{\alpha}{24},\;\;\vspace{0.3 cm}\\
\displaystyle b_2=\left(\frac{11}{2880}+\frac{\alpha}{1152}\right)\alpha,\vspace{0.3 cm}\\
\displaystyle
b_3=-\left(\frac{191}{362880}+\frac{11\alpha}{69120}+\frac{\alpha^2}{82944}\right)\alpha,\,\,\vspace{0.3 cm}\\
\displaystyle
b_4=\left(\frac{2497}{29030400}+\frac{10181\alpha}{348364800}+\frac{11\alpha^2}{3317760}
+\frac{\alpha^3}{7962624}
\right)\alpha.
\end{array}$$

{\it {\bf Remark 1:}
It follows from the proofs in \cite{HSC,MT,TG}, one can see that the conditions stated in Theorem
2 can be weakened as $u(x,t)\in
\mathscr{C}^{2n+\alpha}(\mathrm{R})$ with respect to $x$, where
$$\begin{array}{lll}
\displaystyle
\mathscr{C}^{2n+\alpha}(\mathrm{R})=\left\{u|u\in L_1(\mathrm{R}),
\;and\;\int_{\mathrm{R}}\left(1+|\omega|\right)^{2n+\alpha}|\widehat{u}(\omega,t)|\mathrm{d}\omega<\infty\right\}.
\end{array}
$$}

{\it{\bf Remark 2:} It's worth noting that some suitable smoothness conditions with respect to $x$ for a given
function $u(x,t)$ are necessary and can't be dropped. If one does that, the expected convergence order (accuracy)
in (4) can not be achieved, the fact will be further proved in the Section 4.}

{\it {\bf Remark 3:} It is quite clear that equ. (4) holds for a function $u(x,t)$ defined on
$\mathrm{R}$ with fixed $t\in[0,T]$. If function $u(x,t)$ is defined on the bounded interval $[l,L]$
 with respect to $x$ and satisfies $u(l,t)=u(L,t)=0$, then $u(x,t)$ can be zero-extend smooth outside of the domain
 boundary by defining
$$
\widetilde{u}(x,t)=\left\{
\begin{array}{lll}\displaystyle
 u(x,t),\;(x,t)\in[l,L]\times[0,T];\vspace{0.2 cm}\\
0,\;otherwise.
\end{array}
\right.
$$
In this case, function $\widetilde{u}(x,t)$ is defined on $\mathrm{R}$. Let us suppose that
the extended function $\widetilde{u}(x,t)$ satisfy the conditions of The Theorem 2 or Remark 1,
 then the high-order fractional-compact numerical approximation formula (4) can be rewritten as the following form
for any $x\in[l,L]$,
$$\begin{array}{lll}
\displaystyle \frac{\partial^\alpha
u(x,t)}{\partial{|x|^\alpha}}&=&\displaystyle\left(\delta_x^{0}-b_{n-1} \delta_x^{2n-2}\right)^{-1}\left(-\frac{1}{h^{\alpha}}\sum\limits_{\ell=0}^
{n-2}b_{\ell}\,\delta_x^{2\ell}\sum\limits_{k=-\frac{L-x}{h}+\ell}^
{\frac{x-l}{h}-\ell}g_k^{(\alpha)}u(x-kh,t)\right)\vspace{0.2 cm}\\&&\displaystyle+\mathcal{\mathcal{O}}(h^{2n}),\;\;n\in \textbf{N}^{+}.
\end{array}
$$}

{\it {\bf Remark 4:}  If $\alpha\rightarrow2$, equs. (11)--(14) are reduced  into the following (classical)
 compact formulas,
$$\begin{array}{lll}
\displaystyle \frac{\partial^2
u(x,t)}{\partial{x^2}}=\frac{1}{h^2}\left({\delta_x^{0}+\frac{1}{12}\delta_x^2}\right)^{-1}\delta_x^2u(x,t)
+{\mathcal{O}}(h^4),
\end{array}
$$
$$\begin{array}{lll}
\displaystyle \frac{\partial^2
u(x,t)}{\partial{x^2}}=\frac{1}{h^2}\left(\delta_x^{0}-\frac{1}{90}\delta_x^4\right)^{-1}\left(\delta_x^{0}-\frac{1}{12}\delta_x^2\right)\delta_x^2u(x,t)
+ {\mathcal{O}}(h^6),
\end{array}
$$
$$\begin{array}{lll}
\displaystyle \frac{\partial^2
u(x,t)}{\partial{x^2}}=\frac{1}{h^2}\left(\delta_x^{0}+\frac{1}{560}\delta_x^6\right)^{-1}\left(\delta_x^{0}-\frac{1}{12}\delta_x^2+\frac{1}{90}\delta_x^4\right)
\delta_x^2u(x,t)
+  {\mathcal{O}}(h^8),
\end{array}
$$
$$\begin{array}{lll}
\displaystyle \frac{\partial^2
u(x,t)}{\partial{x^2}}=\frac{1}{h^2}\left(\delta_x^{0}-\frac{1}{3150}\delta_x^8\right)^{-1}\left(\delta_x^{0}-\frac{1}{12}\delta_x^2+\frac{1}{90}\delta_x^4-\frac{1}{560}\delta_x^6\right)
\delta_x^2u(x,t) +  {\mathcal{O}}(h^{10}).
\end{array}
$$}

\section{The application of 4th-order fractional-compact numerical differential formula
to the Riesz spatial telegraph equation}

\quad In this section, we study the Riesz spatial telegraph equation in the following form
 $$
\begin{array}{rrr} \displaystyle
\displaystyle \frac{\partial^2 u(x,t)}{\partial
t^2}+\nu\frac{\partial u(x,t)}{\partial t}=\kappa^2
\frac{\partial^\alpha u(x,t)}{\partial{|x|^\alpha}}
  +f(x,t), \displaystyle\,\,\,l<x<L,\,\,\,0< t\leq T,$$
\end{array}\eqno(15)
$$
with the initial and
boundary value conditions below
$$
\begin{array}{lll}\displaystyle
u(x,0)=\varphi(x),\,\,\,l\leq x \leq L,
\end{array}
$$
$$
\begin{array}{lll}\displaystyle
\frac{\partial u(x,0)}{\partial t}=\phi(x),\,\,\,l\leq x \leq L,
\end{array}
$$
$$
\begin{array}{lll}\displaystyle
\displaystyle u(l,t)=0,\,\,\,0\leq t
\leq T,
\end{array}
$$
$$
\begin{array}{lll}\displaystyle
\displaystyle u(L,t)=0,\,\,\,0\leq t
\leq T,
\end{array}
$$
where
$1<\alpha<2$, $\nu>0$ and $\kappa^2$ are two constants.
It is well known that the wave phenomena and propagation of electric signals
in a cable transmission
line commonly is often reflected by the classical telegraph equations. Regretfully,
the classical telegraph equations cannot describe adequately
the abnormal diffusion phenomena during the finite long transmits progress, where
the voltage or current waves possibly occur \cite{A,BZM}. Due to the memory and
 hereditary properties of different substances can be described by using
the fractional derivatives
 and integrals, so, it is necessary to study
the fractional telegraph equations.  In fact, the fractional telegraph equations
 are the mixture models between diffusion and wave
propagation \cite{J,LGY}, hence, they are more suitable for characterizing transmission and propagation of electrical signals
than the classical ones \cite{OB,OZ}.

In order to guarantee the convergence order of the finite difference scheme
$\mathcal{O}(\tau^4+h^4)$ by using the 4th-order compact and fractional-compact (taking $n=2$ in equ (4)) numerical approximation
formulas for temporal derivative and spatial Riesz derivative, respectively.
Here, we simply assume that equ. (15) has a unique solution $u(x,t)\in C_{x,t}^{6,5}(\Omega)(\Omega=[l,L]\times[0,T])$
subject to the above initial and
boundary value conditions \cite{CJLT,D1,D2}, in which
$$
\begin{array}{lll}\displaystyle
\displaystyle C_{x,t}^{6,5}(\Omega)=\left\{u(x,t)\left|\right.\frac{\partial^p u(x,t)}
{\partial x^p},\frac{\partial^q u(x,t)}{\partial t^q}\in C(\Omega),0\leq p\leq6,0\leq q\leq 5\right\}.
\end{array}
$$

\subsection{Construction of the numerical scheme}
\quad
Throughout the paper,
 let $x_j = l+jh$ in $[l,L]$, $j= 0, 1,
\ldots, M$, and $t_s = s\tau$ in $[0,T]$, $s= 0, 1, \ldots, N$, where $h =
(L-l)/M$ and $\tau = T/N$, $M$ and $N$ are two positive integers.

For convenience, denote
$$
\begin{array}{lll} \displaystyle
\displaystyle
\mu_tu(x,t)=u\left(x,t+\frac{\tau}{2}\right)+u\left(x,t-\frac{\tau}{2}\right),\;\;
\delta_t u(x,t)=u\left(x,t+\frac{\tau}{2}\right)-u\left(x,t-\frac{\tau}{2}\right),\vspace{0.2 cm}\\
\mu_t\delta_t u(x,t)=u(x,t+\tau)-u(x,t-\tau),\;\;
\delta_t^2u(x,t)=u(x,t+\tau)-2u(x,t)+u(x,t-\tau),\vspace{0.2 cm}\\
\delta_x^{\alpha}u(x_j,t)=\sum\limits_{k=-(M-j)}^
{j}g_k^{(\alpha)}u(x-kh,t).
\end{array}
$$

Next, we consider the equation
 $$
\begin{array}{rrr} \displaystyle
\displaystyle \frac{\partial^2 u(x,t)}{\partial
t^2}+\nu\frac{\partial u(x,t)}{\partial t}=g(x,t).
\end{array}
$$
From \cite{TLY}, one can obtain a 4th-order compact
 difference scheme for the above equation,
$$
\begin{array}{rrr} \displaystyle
\displaystyle H_t^{-1}J_tu(x,t)=g(x,t)+\mathcal{O}(\tau^4),
\end{array}\eqno(16)
$$
where $H_t$ and $J_t$
are two operators and defined by
$$
\begin{array}{rrr} \displaystyle
\displaystyle H_t=1+\frac{1}{12}\left(\delta_t^2+\frac{\nu \tau}{2}\mu_t\delta_t\right),\;\;
 J_t=\left(\frac{\nu^2}{12}+\frac{1}{\tau^2}\right)\delta_t^2+\frac{\nu }{2\tau}\mu_t\delta_t.
\end{array}
$$

Replacing $g(x,t)$ by function $\kappa^2
\frac{\partial^\alpha u(x,t)}{\partial{|x|^\alpha}}+f(x,t)$ in equ. (16) gives
$$
\begin{array}{rrr} \displaystyle
\displaystyle H_t^{-1}J_tu(x,t)=\kappa^2
\frac{\partial^\alpha u(x,t)}{\partial{|x|^\alpha}}
  +f(x,t)+\mathcal{O}(\tau^4).
\end{array}\eqno(17)
$$

For Riesz derivative $
\frac{\partial^\alpha u(x,t)}{\partial{|x|^\alpha}}$, we choose the 4th-order fractional-compact numerical differential formula
(11). It follows from (17) that
$$
\begin{array}{rrr} \displaystyle
\displaystyle H_t^{-1}J_tu(x,t)=\kappa^2
\left({\delta_x^0+\frac{\alpha}{24}\delta_x^2}\right)^{-1}\left(-\frac{\Delta_h^{\alpha}u(x,t)}{h^{\alpha}}\right)
  +f(x,t)+R(x,t),
\end{array}\eqno(18)
$$
where there exists a constant $c_2$ such that
 $R(x,t)\leq c_2 (\tau^4+h^4)$.

Neglecting the high order term $R(x,t)$ and letting $u_j^s$
be the approximate solution of function $u(x_j,t_s),$ then one has
$$
\begin{array}{rrr} \displaystyle
\displaystyle H_t^{-1}J_tu_j^s=-\frac{\kappa^2}{h^\alpha}
\left({\delta_x^0+\frac{\alpha}{24}\delta_x^2}\right)^{-1}\delta_x^{\alpha}u_{j}^s
  +f_j^s,\vspace{0.2 cm}\\\;\;j=1,\ldots,M-1,\;s=1,2,\ldots,N-1.
\end{array}\eqno(19)
$$

As for the finite difference scheme (19), it is a three-level one which
needs knowing the approximate value at $t = t_1$. Here, we provide a way to compute function value $u(x,t)$
at the first time level. For convenience, firstly denote $\frac{\partial^m u(x,t)}{\partial t^m}|_{t=0}=\frac{\partial^m u(x,0)}{\partial t^m}$, $m=0,1,\ldots$.

From equ. (15), one has
$$\begin{array}{lll}
\displaystyle
\frac{\partial^2 u(x,0)}{\partial
t^2}&=&\displaystyle\kappa^2
\frac{\partial^\alpha u(x,0)}{\partial{|x|^\alpha}}-\nu\frac{\partial u(x,0)}{\partial t}
  +f(x,0)\vspace{0.2 cm}\\&=&\displaystyle\kappa^2
\frac{\partial^\alpha \varphi(x)}{\partial{|x|^\alpha}}-\nu\phi(x)
  +f(x,0),
\end{array}
$$
and
$$\begin{array}{lll}
\displaystyle
\frac{\partial^3 u(x,0)}{\partial
t^3}&=&\displaystyle\kappa^2\frac{\partial}{\partial t}\left(
\frac{\partial^\alpha u(x,t)}{\partial{|x|^\alpha}}\right)\hspace{-0.2 cm}\bigg|_{t=0}-\nu\frac{\partial^2 u(x,0)}{\partial t^2}
  +\frac{\partial f(x,0)}{\partial t}\vspace{0.2 cm}\\&=&\displaystyle
  \frac{\kappa^2}{\tau}\left(\frac{\partial^\alpha u(x,t_1)}{\partial{|x|^\alpha}}-\frac{\partial^\alpha \varphi(x)}{\partial{|x|^\alpha}}\right)
-\nu\left(\kappa^2
\frac{\partial^\alpha \varphi(x)}{\partial{|x|^\alpha}}-\nu\phi(x)
  +f(x,0)\right)
 \vspace{0.2 cm}\\&&\displaystyle+\frac{\partial f(x,0)}{\partial t}+\mathcal{O}(\tau).
\end{array}
$$

Next, we obtain the following equation
by substituting the above two equations into the Taylor's series expansion of function $u(x,t_1)$ at point $t=0$,
$$\begin{array}{lll}
\displaystyle
u(x,t_1)&=&\displaystyle u(x,0)+\tau\frac{\partial u(x,0)}{\partial
t}+\frac{\tau^2}{2}
\frac{\partial^2 u(x,0)}{\partial
t^2}+\frac{\tau^3}{6}
\frac{\partial^3 u(x,0)}{\partial
t^3}+\mathcal{O}(\tau^4)\vspace{0.2 cm}\\&=&\displaystyle
\varphi(x)+\frac{\tau}{6}(6-3\nu\tau+\nu^2\tau^2)\phi(x)+\frac{\tau^2}{6}(3-\nu\tau)f(x,0)
+\frac{\tau^3}{6}\frac{\partial f(x,0)}{\partial t}\vspace{0.2 cm}\\&&\displaystyle+
\frac{\kappa^2\tau^2}{6}(2-\nu^2)\left({\delta_x^0+\frac{\alpha}{24}\delta_x^2}\right)^{-1}\left(-\frac{\Delta_h^{\alpha}\varphi(x)}{h^{\alpha}}\right)
\vspace{0.2 cm}\\&&\displaystyle+
\frac{\kappa^2\tau^2}{6}\left({\delta_x^0+\frac{\alpha}{24}\delta_x^2}\right)^{-1}\left(-\frac{\Delta_h^{\alpha}u(x,t_1)}{h^{\alpha}}\right)
+R(x,t_1),
\end{array}
$$
where there exists a constant $c_1$ such that $R(x,t_1)\leq c_1(\tau^4+h^4)$.

So, the first level value $u(x_j,t_1)$ can be numerical determined by the
equation
$$\begin{array}{lll}
&&\displaystyle
\left({\delta_x^0+\frac{\alpha}{24}\delta_x^2}\right)u_j^1+
\frac{\kappa^2\tau^2}{6h^\alpha}\delta_x^{\alpha}u_{j}^1
\vspace{0.2 cm}
\\&=&\displaystyle\left({\delta_x^0+\frac{\alpha}{24}\delta_x^2}\right)\left(
\varphi(x_j)+\frac{\tau}{6}(6-3\nu\tau+\nu^2\tau^2)\phi(x_j)+\frac{\tau^2}{6}(3-\nu\tau)f(x_j,0)\right.\vspace{0.2 cm}
\\ &&\displaystyle\left.
+\frac{\tau^3}{6}\frac{\partial f(x_j,0)}{\partial t}\right)-
\frac{\kappa^2\tau^2}{6h^\alpha}(2-\nu^2)
\sum\limits_{k=-(M-j)}^
{j}g_k^{(\alpha)}\varphi(x_{j-k}),\;\;j=1,2,\ldots, M-1.
\end{array}
$$

Therefore, we obtain a high-order finite difference scheme for equ. (15) as follows,
$$\left\{
\begin{array}{lll} \displaystyle
\displaystyle \left({\delta_x^0+\frac{\alpha}{24}\delta_x^2}\right)J_tu_j^s=-\frac{\kappa^2}{h^\alpha}
H_t\delta_x^{\alpha}u_{j}^s
  +\left({\delta_x^0+\frac{\alpha}{24}\delta_x^2}\right)H_tf_j^s,\;\;j=1,\ldots,M-1,\vspace{0.2 cm}\\\;\hspace{8.5 cm}s=2,\ldots,N-1,\vspace{0.2 cm}\\
  \displaystyle
  u_j^0=\varphi(x_j),
\;\;j=1,\ldots,M-1,\vspace{0.2 cm}\\
   \displaystyle
\left({\delta_x^0+\frac{\alpha}{24}\delta_x^2}\right)u_j^1+
\frac{\kappa^2\tau^2}{6h^\alpha}\delta_x^{\alpha}u_{j}^1=\displaystyle\left({\delta_x^0+\frac{\alpha}{24}\delta_x^2}\right)\left(
\varphi(x_j)+\frac{\tau}{6}(6-3\nu\tau+\nu^2\tau^2)\phi(x_j)\right.\vspace{0.2 cm}
\\ \displaystyle\left.+\frac{\tau^2}{6}(3-\nu\tau)f(x_j,0)
+\frac{\tau^3}{6}\frac{\partial f(x_j,0)}{\partial t}\right)-
\frac{\kappa^2\tau^2}{6h^\alpha}(2-\nu^2)
\sum\limits_{k=-(M-j)}^
{j}g_k^{(\alpha)}\varphi(x_{j-k}),\vspace{0.2 cm}\\\;\hspace{8.5 cm}j=1,2,\ldots, M-1.
\vspace{0.2 cm}\\
\displaystyle
u_0^s=u_M^s=0,
\;
\;s=0,1,\ldots,N.
\end{array}\right.\eqno(20)
$$

\subsection{Stability analysis}

Denote
$$\begin{array}{lll} \displaystyle
V_h=\{v|v=\{v_j\},\;j=0,1,\ldots,M\},\;\; \mathring{V}_h=\{v|v\in S_h,v_0=v_M=0\},
\end{array}
$$
then for any $u,v\in\mathring{S}_h$, we can define the following inner products:
$$\begin{array}{lll} \displaystyle
(u,v)=h\sum_{j=1}^{M-1}u_jv_j,\;\;(\delta_xu,\delta_xv)=h\sum_{j=1}^{M-1}(\delta_xu_{j-\frac{1}{2}})(\delta_xv_{j-\frac{1}{2}})
\end{array}
$$
and the following norms:
$$\begin{array}{lll} \displaystyle
||u||=\sqrt{(u,u)},\;\;\;\;
||\delta_xu||=
\sqrt{
(\delta_xu,\delta_xv)}.
\end{array}
$$
Meanwhile, we also easily obtain
$$\begin{array}{lll} \displaystyle
(\delta_x^2u,v)=-(\delta_xu,\delta_xv),
\end{array}
$$
and the inverse estimate formula
$$\begin{array}{lll} \displaystyle
||\delta_xu||^2\leq4||u||^2.
\end{array}
$$

Next, we give several lemmas as follows:

{\it{\bf Lemma 2.} (Gerschgorin theorem \cite{QSS})  Let $Z=(z_{js})$ be a complex matrix of order $M-1$, and
$$
\begin{array}{lll}
\displaystyle\displaystyle Y_j=\sum_{s \ne j} {|z_{js}|},\quad j=1,
2,\dots, M-1.
\end{array}
$$
Let $X_j$ be the closed disc centered at $z_{jj}$ with radius
$Y_j$:
$$
\begin{array}{lll}
\displaystyle X_j=\{\, x \in \mathbf C: |x-z_{jj}|\le Y_j\,\}.
\end{array}
$$
Then all the eigenvalues of matrix $Z$ belong to
$\bigcup\limits_{j = 1}^{M - 1}X_j $.}

{\it{\bf Lemma 3.} The operator $\left({\delta_x^0+\frac{\alpha}{24}\delta_x^2}\right)$ is self-adjoint, i.e., for any $u,v\in \mathring{V}_h$, one has
$$
\begin{array}{lll}
\displaystyle
\left(\left({\delta_x^0+\frac{\alpha}{24}\delta_x^2}\right)u,v\right)=\left(u,\left({\delta_x^0+\frac{\alpha}{24}\delta_x^2}\right)v\right).
\end{array}
$$
}

{\bf Proof.} Note that
$$
\begin{array}{lll}
\displaystyle
\left(\left({\delta_x^0+\frac{\alpha}{24}\delta_x^2}\right)u,v\right)&=&\displaystyle\left(u,v\right)+\frac{\alpha}{24}\left(\delta_x^2u,v\right)=
\left(u,v\right)-\frac{\alpha}{24}\left(\delta_xu,\delta_xv\right)\vspace{0.2 cm}\\
&=&\displaystyle\left(u,v\right)+\frac{\alpha}{24}\left( u,\delta_x^2v\right)
=\left(u,\left({\delta_x^0+\frac{\alpha}{24}\delta_x^2}\right)v\right),
\end{array}
$$
this finishes the proof.\QEDA

{\it{\bf Lemma 4.} For any $v\in \mathring{V}_h$, one has the following estimate:
$$
\begin{array}{lll}
\displaystyle
\left(1-\frac{\alpha}{6}\right)||v||^2\leq\left(\left({\delta_x^0+\frac{\alpha}{24}\delta_x^2}\right)v,v\right)\leq ||v||^2.
\end{array}
$$
}

{\bf Proof.} On the one hand
$$
\begin{array}{lll}
\displaystyle
\left(\left({\delta_x^0+\frac{\alpha}{24}\delta_x^2}\right)u,v\right)=
||v||^2-\frac{\alpha}{24}||\delta_xv||^2\leq||v||^2,
\end{array}
$$
on the other hand, we also have
$$
\begin{array}{lll}
\displaystyle
\left(\left({\delta_x^0+\frac{\alpha}{24}\delta_x^2}\right)u,v\right)=
||v||^2-\frac{\alpha}{24}||\delta_xv||^2\geq||v||^2-\frac{\alpha}{6}||v||^2,
\end{array}
$$
thus, the proof was completed.\QEDA

{\it{\bf Lemma 5.} The fractional operator $\delta_x^\alpha$ is self-adjoint, i.e., for any $u,v\in \mathring{V}_h$, one has
$$
\begin{array}{lll}
\displaystyle
\left(\delta_x^\alpha u,v\right)=\left(u,\delta_x^\alpha v\right).
\end{array}
$$
}

{\bf Proof.} For convenience, we define the following
 matrix $A_\alpha$ as the corresponding matrix associated with the fractional operator $\delta_x^\alpha$:
$$ \displaystyle A_\alpha= \left(
  \begin{array}{ccccccc}
    g_0^{(\alpha)} & g_{-1}^{(\alpha)} & g_{-2}^{(\alpha)} & \cdots & g_{4-M}^{(\alpha)} & g_{3-M}^{(\alpha)} & g_{2-M}^{(\alpha)} \vspace{0.2 cm}\\
    g_{1}^{(\alpha)} & g_{0}^{(\alpha)}& g_{-1}^{(\alpha)}& g_{-2}^{(\alpha)} & \cdots &g_{4-M}^{(\alpha)} &    g_{3-M}^{(\alpha)} \vspace{0.2 cm}\\
    g_{2}^{(\alpha)} & g_{1}^{(\alpha)} & g_{0}^{(\alpha)} & g_{-1}^{(\alpha)} &g_{-2}^{(\alpha)} & \cdots &    g_{4-M}^{(\alpha)} \vspace{0.2 cm}\\
    \vdots & \ddots & \ddots & \ddots & \ddots & \ddots & \vdots\vspace{0.2 cm}\\
    g_{M-4}^{(\alpha)} & \cdots & g_{2}^{(\alpha)} & g_{1}^{(\alpha)} & g_{0}^{(\alpha)} & g_{-1}^{(\alpha)} & g_{-2}^{(\alpha)}\vspace{0.2 cm}\\
    g_{M-3}^{(\alpha)} & g_{M-4}^{(\alpha)} & \cdots & g_{2}^{(\alpha)} & g_{1}^{(\alpha)} & g_{0}^{(\alpha)} & g_{-1}^{(\alpha)} \vspace{0.2 cm}\\
    g_{M-2}^{(\alpha)}& g_{M-3}^{(\alpha)} & g_{M-4}^{(\alpha)} & \cdots & g_{2}^{(\alpha)} & g_{1}^{(\alpha)} & g_{0}^{(\alpha)} \vspace{0.2 cm}\\
  \end{array}
\right).
$$
The matrix form of $\delta_x^\alpha v_j$ $(j=1,2,\ldots, M-1)$ can be denoted by $A_\alpha v$.
The direct calculation leads to
$$
\begin{array}{lll}
\displaystyle
\left(\delta_x^\alpha u,v\right)=v^{T}A_{\alpha}u=\left(v^{T}A_{\alpha}u\right)^{T}=u^{T}A_{\alpha}v=
\left(u,\delta_x^\alpha v\right).
\end{array}
$$
This ends the proof.\QEDA

{\it{\bf Lemma 6.} The fractional operator $\delta_x^\alpha$ is positive definite, i.e., for any $v\in \mathring{V}_h$, one has
$$
\begin{array}{lll}
\displaystyle
\left(\delta_x^\alpha v,v\right)>0.
\end{array}
$$
Furthermore, there also exists a positive definite operator $\widetilde{\delta}_x^\alpha$, such that
$$
\begin{array}{lll}
\displaystyle
\left(\delta_x^\alpha v,v\right)=\left(\widetilde{\delta}_x^\alpha v,\widetilde{\delta}_x^\alpha v\right).
\end{array}
$$
}

{\bf Proof.} In view of the properties of the coefficients $g_k^{(\alpha)}$, we easily know that matrix $A_\alpha$ is symmetric.

Next, we prove the definite positivity of $A_\alpha$. From Lemma 2 and Theorem 2, we know that the eigenvalues $\lambda_j(A_\alpha)\;(j=1,2,\ldots,M-1)$
of the matrix $A_\alpha$ satisfy,
$$
\begin{array}{lll}
\displaystyle
|\lambda_j(A_\alpha)-g_0^{(\alpha)}|\leq\sum\limits_{k=-(M-2),k\neq0}^
{M-2}|g_k^{(\alpha)}|,
\end{array}
$$
i.e.,
$$
\begin{array}{lll}
\displaystyle
-\sum\limits_{k=-(M-2),k\neq0}^
{M-2}|g_k^{(\alpha)}|&\leq&\displaystyle\lambda_j(A_\alpha)-g_0^{(\alpha)}\leq\sum\limits_{k=-(M-2),k\neq0}^
{M-2}|g_k^{(\alpha)}|\vspace{0.2 cm}\\
&<&\displaystyle\sum\limits_{k=-\infty,k\neq0}^
{\infty}|g_k^{(\alpha)}|\leq\frac{2^{1+\alpha}}{\pi},
\end{array}
$$
or,
$$
\begin{array}{lll}
\displaystyle
 g_0^{(\alpha)}-\sum\limits_{k=-(M-2),k\neq0}^
{M-2}|g_k^{(\alpha)}|\leq\lambda_j(A_\alpha)< g_0^{(\alpha)}+\frac{2^{1+\alpha}}{\pi}.
\end{array}
$$
Noting that
$$
\begin{array}{lll}
\displaystyle
 g_0^{(\alpha)}-\sum\limits_{k=-(M-2),k\neq0}^
{M-2}|g_k^{(\alpha)}|=\sum\limits_{k=-(M-2)}^
{M-2}g_k^{(\alpha)}>\sum\limits_{k=-\infty}^
{\infty}g_k^{(\alpha)}=0,
\end{array}
$$
then one has
$$
\begin{array}{lll}
\displaystyle
0<\lambda_j(A_\alpha)< g_0^{(\alpha)}+\frac{2^{1+\alpha}}{\pi},
\end{array}
$$
that is, the matrix $A_\alpha$ is positive definite, i.e.,
$$
\begin{array}{lll}
\displaystyle
\left(\delta_x^\alpha v,v\right)=v^{T}A_\alpha v>0.
\end{array}
$$

Because of the matrix $A_\alpha$ is symmetric and positive definite, then there must be exists a
symmetric and positive definite matrix $B_\alpha$, such that
$$
\begin{array}{lll}
\displaystyle
A_\alpha =B_\alpha^2.
\end{array}
$$

Thus,
$$
\begin{array}{lll}
\displaystyle
\left(\delta_x^\alpha v,v\right)=v^{T}A_\alpha v=(v^{T}B_\alpha)(B_\alpha v)=(B_\alpha v)^{T}(B_\alpha v)
=\left(\widetilde{\delta}_x^\alpha v,\widetilde{\delta}_x^\alpha v\right),
\end{array}
$$
where $\widetilde{\delta}_x^\alpha$ is the associated fractional difference operator of matrix $B_\alpha$. This completes the proof.\QEDA

{\it{\bf Lemma 7.} For any $v\in \mathring{V}_h$, one has the following estimate:
$$
\begin{array}{lll}
\displaystyle
 \lambda_{\min}(A_\alpha)||v||^2\leq\left(\delta_x^\alpha v,v\right)\leq \lambda_{\max}(A_\alpha)||v||^2,
\end{array}
$$
where $\lambda_{\min}(A_\alpha)$ and $\lambda_{\max}(A_\alpha)$ denote the
 smallest and
 the largest eigenvalues of matrix $A_\alpha$, respectively.}

{\bf Proof.}
For any symmetric matrix $A_\alpha$, we can
use the property of Rayleigh-Ritz ratio \cite{HJ}, i.e., for any vector $v\neq0$, we have
$$
\begin{array}{lll}
\displaystyle
 \lambda_{\min}(A_\alpha)\leq\frac{\left(A_\alpha v,v\right)}{(v,v)}\leq \lambda_{\max}(A_\alpha).
\end{array}
$$

Because of $A_\alpha v$ denote the matrix form of $\delta_x^\alpha v_j$ $(j=1,2,\ldots, M-1)$,
then one has
$$
\begin{array}{lll}
\displaystyle
 \lambda_{\min}(A_\alpha)||v||^2\leq\left(\delta_x^\alpha v,v\right)\leq \lambda_{\max}(A_\alpha)||v||^2.
\end{array}
$$
This finishes the proof.\QEDA

{\it{\bf Lemma 8.} If
$$
\begin{array}{lll} \displaystyle
\frac{\tau}{h^\alpha}\left(\frac{\tau\nu^2}{12}+\frac{1}{\tau}\right)^{-1}\leq\frac{6-\alpha}{2\kappa^2g_0^{(\alpha)}},
\end{array}\eqno(21)
$$
then the operator $\left(\left({\delta_x^0+\frac{\alpha}{24}\delta_x^2}\right)-\frac{\tau\kappa^2}{6h^\alpha}
\left(\frac{\tau\nu^2}{12}+\frac{1}{\tau}\right)^{-1}\delta_x^\alpha\right)$
is positive definite, i.e.,
for any $v\in \mathring{V}_h$, there holds
$$
\begin{array}{lll}
\displaystyle
\left(\left(\left({\delta_x^0+\frac{\alpha}{24}\delta_x^2}\right)-\frac{\tau\kappa^2}
{6h^\alpha}\left(\frac{\tau\nu^2}{12}+\frac{1}{\tau}\right)^{-1}\delta_x^\alpha\right) v,v\right)>0.
\end{array}
$$
}

{\bf Proof.} Denote the incidence matrix of the operator $\left({\delta_x^0+\frac{\alpha}{24}\delta_x^2}\right)$ is
$$ \displaystyle
D_\alpha= \left(
  \begin{array}{cccccc}
   1-\frac{\alpha}{12} & \frac{\alpha}{24}&  0& & \cdots& 0 \vspace{0.2 cm}\\
   \frac{\alpha}{24} &  1-\frac{\alpha}{12} &\frac{\alpha}{24}& 0& &\vdots \vspace{0.2 cm}\\
     & \frac{\alpha}{24} & 1-\frac{\alpha}{12}  & \frac{\alpha}{24} & \ddots&  \vspace{0.2 cm}\\
    \vdots & \ddots & \ddots & \ddots & \ddots & \vdots \vspace{0.4 cm}\\
    0& \ldots&  &&   \frac{\alpha}{24}  & 1-\frac{\alpha}{12}\vspace{0.2 cm}\\
  \end{array}
\right).
$$

Obviously, the incidence matrix of the operator
$$\left(\left({\delta_x^0+\frac{\alpha}{24}\delta_x^2}\right)-\frac{\tau\kappa^2}
{6h^\alpha}\left(\frac{\tau\nu^2}{12}+\frac{1}{\tau}\right)^{-1}\delta_x^\alpha\right)$$ is
$$
\begin{array}{lll}
\displaystyle
G_\alpha=D_\alpha-\frac{\tau\kappa^2}{6h^\alpha}\left(\frac{\tau\nu^2}{12}+\frac{1}{\tau}\right)^{-1} A_\alpha,
\end{array}
$$
and its eigenvalues satisfy:
$$
\begin{array}{lll}
\displaystyle
\left|\lambda_j(G_\alpha)-\left[\left(1-\frac{\alpha}{12}\right)-\frac{\tau\kappa^2}{6h^\alpha}\left(\frac{\tau\nu^2}{12}+\frac{1}{\tau}\right)^{-1} g_0^{(\alpha)}\right]\right|
\vspace{0.2 cm}\\\displaystyle
\leq
\frac{\alpha}{12}+\frac{\tau\kappa^2}{6h^\alpha}\left(\frac{\tau\nu^2}{12}+\frac{1}{\tau}\right)^{-1}\sum\limits_{k=-(M-2),k\neq0}^
{M-2}|g_k^{(\alpha)}|
\vspace{0.2 cm}\\\displaystyle
<
\frac{\alpha}{12}+\frac{\tau\kappa^2}{6h^\alpha}\left(\frac{\tau\nu^2}{12}+\frac{1}{\tau}\right)^{-1}\sum\limits_{k=-\infty,k\neq0}^
{\infty}|g_k^{(\alpha)}|,
\end{array}
$$
so,
$$
\begin{array}{lll}
\displaystyle
\lambda_j(G_\alpha)&>&\displaystyle
\left(1-\frac{\alpha}{6}\right)-\frac{\tau\kappa^2}{6h^\alpha}\left(\frac{\tau\nu^2}{12}+\frac{1}{\tau}\right)^{-1}\left(g_0^{(\alpha)}+\sum\limits_{k=-\infty,k\neq0}^
{\infty}|g_k^{(\alpha)}|\right)\geq0
\end{array}
$$
under the condition (21). In other words,
$$
\begin{array}{lll}
\displaystyle
\left(\left(\left({\delta_x^0+\frac{\alpha}{24}\delta_x^2}\right)-\frac{\tau\kappa^2}
{6h^\alpha}\left(\frac{\tau\nu^2}{12}+\frac{1}{\tau}\right)^{-1}\delta_x^\alpha\right) v,v\right)=v^{T}G_\alpha
v>0.
\end{array}
$$
This finishes the proof.
\QEDA

{\it{\bf Lemma 9.} For any $v\in \mathring{V}_h$, there holds:
$$
\begin{array}{lll}
\displaystyle
\left\|\left({\delta_x^0+\frac{\alpha}{24}\delta_x^2}\right)v\right\|\leq ||v||.
\end{array}
$$
}

{\bf Proof.} According to the definition, one easily obtain

$$
\begin{array}{lll}
\displaystyle\left\|\left({\delta_x^0+\frac{\alpha}{24}\delta_x^2}\right)v\right\|^2&=&\displaystyle
\left(\left({\delta_x^0+\frac{\alpha}{24}\delta_x^2}\right)v,\left({\delta_x^0+\frac{\alpha}{24}\delta_x^2}\right)v\right)
\vspace{0.2 cm}\\\displaystyle
&=&\displaystyle
||v||^2-\frac{\alpha}{12}||\delta_xv||^2+\frac{\alpha^2}{576}||\delta_x^2v||^2

\vspace{0.2 cm}\\\displaystyle
&\leq&\displaystyle
||v||^2+\frac{\alpha^2-12\alpha}{144}||\delta_xv||^2
\vspace{0.2 cm}\\\displaystyle
&\leq&\displaystyle
||v||^2.
\end{array}
$$
This completes the proof.
\QEDA

Now, we give the stability analysis in details.

 {\it{\bf Theorem 3.}The difference scheme  (20) is stable to the initial value under condition (21).}

{\bf Proof.}
Suppose that $v_j^s$ is the solution of the following difference equation:
$$\left\{
\begin{array}{lll} \displaystyle
\displaystyle \left({\delta_x^0+\frac{\alpha}{24}\delta_x^2}\right)J_tv_j^s=-\frac{\kappa^2}{h^\alpha}
H_t\delta_x^{\alpha}v_{j}^s
  +\left({\delta_x^0+\frac{\alpha}{24}\delta_x^2}\right)H_tf_j^s,\;\;j=1,\ldots,M-1,\vspace{0.2 cm}\\ \hspace{8.5 cm}\;s=2,\ldots,N-1,\vspace{0.2 cm}\\
  \displaystyle
  v_j^0=\varphi(x_j)+\rho_j^0,
\;\;j=1,\ldots,M-1,\vspace{0.2 cm}\\
   \displaystyle
\left({\delta_x^0+\frac{\alpha}{24}\delta_x^2}\right)v_j^1+
\frac{\kappa^2\tau^2}{6h^\alpha}\delta_x^{\alpha}v_{j}^1=\displaystyle\left({\delta_x^0+\frac{\alpha}{24}\delta_x^2}\right)\left(
\varphi(x_j)+\frac{\tau}{6}(6-3\nu\tau+\nu^2\tau^2)\phi(x_j)\right.\vspace{0.2 cm}
\\ \displaystyle\left.+\frac{\tau^2}{6}(3-\nu\tau)f(x_j,0)
+\frac{\tau^3}{6}\frac{\partial f(x_j,0)}{\partial t}\right)-
\frac{\kappa^2\tau^2}{6h^\alpha}(2-\nu^2)
\sum\limits_{k=-(M-j)}^
{j}g_k^{(\alpha)}\varphi(x_{j-k})+\rho_j^1,\vspace{0.2 cm}\\ \hspace{8.5 cm}\;\;j=1,2,\ldots, M-1.
\vspace{0.2 cm}\\
\displaystyle
v_0^s=v_M^s=0,
\;
\;s=0,1,\ldots,N.
\end{array}\right.\eqno(22)
$$

Let $\varepsilon_j^s=v_j^s-u_j^s$. Then one obtains the following perturbation equation by using equs. (20) and (22),
$$\left\{
\begin{array}{lll} \displaystyle
\displaystyle \left({\delta_x^0+\frac{\alpha}{24}\delta_x^2}\right)J_t\varepsilon_j^s=-\frac{\kappa^2}{h^\alpha}
H_t\delta_x^{\alpha}\varepsilon_{j}^s,\;\;j=1,\ldots,M-1,\;s=2,\ldots,N-1,\vspace{0.2 cm}\\
  \displaystyle
  \varepsilon_j^0=\rho_j^0,
\;\;j=1,\ldots,M-1,\vspace{0.2 cm}\\
   \displaystyle
\left({\delta_x^0+\frac{\alpha}{24}\delta_x^2}\right)\varepsilon_j^1+
\frac{\kappa^2\tau^2}{6h^\alpha}\delta_x^{\alpha}\varepsilon_{j}^1=\rho_j^1,\;\;j=1,2,\ldots, M-1.
\vspace{0.2 cm}\\
\displaystyle
\varepsilon_0^s=\varepsilon_M^s=0,
\;
\;s=0,1,\ldots,N.
\end{array}\right.\eqno(23)
$$

Firstly, taking the inner product of the 3rd
equation of (23) with $\varepsilon_j^1$, we have
$$
\begin{array}{lll}
\displaystyle
\left(\left({\delta_x^0+\frac{\alpha}{24}\delta_x^2}\right)\varepsilon^1,\varepsilon^1\right)+
\frac{\kappa^2\tau^2}{6h^\alpha}\left(\delta_x^{\alpha}\varepsilon^1,\varepsilon^1\right)=\left(\rho_j^1,,\varepsilon^1\right).
\end{array}
$$

 Application of Lemmas 4 and 6 to the above equation gives,
$$
\begin{array}{lll}
\displaystyle
||\varepsilon^1||\leq\frac{6}{6-\alpha}||\rho^1||.
\end{array}
$$

Next, taking the inner product of the 1st
equation of (23) with $\mu_t\delta_t\varepsilon_j^s$ yields,
$$
\begin{array}{lll}
\displaystyle
\left(\left({\delta_x^0+\frac{\alpha}{24}\delta_x^2}\right)J_t\varepsilon^s,\mu_t\delta_t\varepsilon^s\right)
+\frac{\kappa^2}{h^\alpha}\left(
H_t\delta_x^{\alpha}\varepsilon^s,\mu_t\delta_t\varepsilon^s\right)=0.
\end{array}\eqno(24)
$$

For the first term of equ. (24), we have the following transform by the help of Lemma 3,
$$
\begin{array}{lll}
\displaystyle
\left(\left({\delta_x^0+\frac{\alpha}{24}\delta_x^2}\right)J_t\varepsilon^s,\mu_t\delta_t\varepsilon^s\right)
&=&\displaystyle\left(\frac{\nu^2}{12}+\frac{1}{\tau^2}\right)\left(\left({\delta_x^0+\frac{\alpha}{24}\delta_x^2}\right)\delta_t^2\varepsilon^s,\mu_t\delta_t\varepsilon^s\right)\vspace{0.2 cm}\\
&&\displaystyle+\frac{\nu}{2\tau}\left(\left({\delta_x^0+\frac{\alpha}{24}\delta_x^2}\right)\mu_t\delta_t\varepsilon^s,\mu_t\delta_t\varepsilon^s\right)\vspace{0.2 cm}\\
&=&\displaystyle\left(\frac{\nu^2}{12}+\frac{1}{\tau^2}\right)\left(\left({\delta_x^0+\frac{\alpha}{24}\delta_x^2}\right)\delta_t\varepsilon^{s+\frac{1}{2}},\delta_t\varepsilon^{s+\frac{1}{2}}\right)\vspace{0.2 cm}\\
&&\displaystyle-\left(\frac{\nu^2}{12}+\frac{1}{\tau^2}\right)\left(\left({\delta_x^0+\frac{\alpha}{24}\delta_x^2}\right)\delta_t\varepsilon^{s-\frac{1}{2}},\delta_t\varepsilon^{s-\frac{1}{2}}\right)\vspace{0.2 cm}\\
&&\displaystyle+\frac{\nu}{2\tau}\left(\left({\delta_x^0+\frac{\alpha}{24}\delta_x^2}\right)\mu_t\delta_t\varepsilon^s,\mu_t\delta_t\varepsilon^s\right).
\end{array}
$$

Similarly, by using Lemma 5 to the 2nd term of equ. (24), one has
$$
\begin{array}{lll}
\displaystyle
\frac{\kappa^2}{h^\alpha}\left(
H_t\delta_x^{\alpha}\varepsilon^s,\mu_t\delta_t\varepsilon^s\right)
&=&\displaystyle\frac{\kappa^2}{h^\alpha}\left(\delta_x^\alpha\varepsilon^s,\mu_t\delta_t\varepsilon^s\right)+\frac{\kappa^2}{12h^\alpha}\left(\delta_x^\alpha\delta_t^2\varepsilon^s,\mu_t\delta_t\varepsilon^s\right)
\vspace{0.2 cm}\\
&&\displaystyle+\frac{\nu\tau\kappa^2}{24h^\alpha}\left(\delta_x^\alpha\mu_t\delta_t\varepsilon^s,\mu_t\delta_t\varepsilon^s\right)\vspace{0.2 cm}\\
&=&\displaystyle\frac{\kappa^2}{h^\alpha}\left(\delta_x^\alpha\varepsilon^s,\varepsilon^{s+1}\right)-\frac{\kappa^2}{h^\alpha}\left(\delta_x^\alpha\varepsilon^s,\varepsilon^{s-1}\right)
+\frac{\kappa^2}{12h^\alpha}\left(\delta_x^\alpha\delta_t\varepsilon^{s+\frac{1}{2}},\delta_t\varepsilon^{s+\frac{1}{2}}\right)
\vspace{0.2 cm}\\
&&\displaystyle-\frac{\kappa^2}{12h^\alpha}\left(\delta_x^\alpha\delta_t\varepsilon^{s-\frac{1}{2}},\delta_t\varepsilon^{s-\frac{1}{2}}\right)
+\frac{\nu\tau\kappa^2}{24h^\alpha}\left(\delta_x^\alpha\mu_t\delta_t\varepsilon^s,\mu_t\delta_t\varepsilon^s\right).
\end{array}
$$

Denote
$$
\begin{array}{lll}
\displaystyle
W^s&=&\displaystyle\left(\frac{\tau\nu^2}{12}+\frac{1}{\tau}\right)\left(\left({\delta_x^0+\frac{\alpha}{24}\delta_x^2}\right)\delta_t\varepsilon^{s+\frac{1}{2}},\delta_t\varepsilon^{s+\frac{1}{2}}\right)+
\frac{\tau\kappa^2}{h^\alpha}\left(\delta_x^\alpha\varepsilon^s,\varepsilon^{s+1}\right)\vspace{0.2 cm}\\
&&\displaystyle+\frac{\tau\kappa^2}{12h^\alpha}\left(\delta_x^\alpha\delta_t\varepsilon^{s+\frac{1}{2}},\delta_t\varepsilon^{s+\frac{1}{2}}\right).
\end{array}
$$
Now equ. (24) can be rewritten as
$$
\begin{array}{lll}
\displaystyle
W^s+\frac{\nu}{2}\left(\left({\delta_x^0+\frac{\alpha}{24}\delta_x^2}\right)\mu_t\delta_t\varepsilon^s,\mu_t\delta_t\varepsilon^s\right)
+\frac{\nu\tau^2\kappa^2}{24h^\alpha}\left(\delta_x^\alpha\mu_t\delta_t\varepsilon^s,\mu_t\delta_t\varepsilon^s\right)=W^{s-1}.
\end{array}
$$
Noticing Lemmas 4 and 6, one gets
$$
\begin{array}{lll}
\displaystyle
W^s\leq W^{s-1}\leq\cdots\leq W^0.
\end{array}
$$

In addition, combining Lemma 4 with Lemma 7, we have the following estimate for $W^s$:
$$
\begin{array}{lll}
\displaystyle
W^s&=&\displaystyle\left(\frac{\tau\nu^2}{12}+\frac{1}{\tau}\right)\left(\left({\delta_x^0+\frac{\alpha}{24}\delta_x^2}\right)\delta_t\varepsilon^{s+\frac{1}{2}},\delta_t\varepsilon^{s+\frac{1}{2}}\right)+
\frac{\tau\kappa^2}{h^\alpha}\left(\delta_x^\alpha\varepsilon^s,\varepsilon^{s+1}\right)\vspace{0.2 cm}\\
&&\displaystyle+\frac{\tau\kappa^2}{12h^\alpha}\left(\delta_x^\alpha\delta_t\varepsilon^{s+\frac{1}{2}},\delta_t\varepsilon^{s+\frac{1}{2}}\right)\vspace{0.2 cm}\\
&\leq&\displaystyle2\left(\frac{\tau\nu^2}{12}+\frac{1}{\tau}\right)\left(||\varepsilon^{s+1}||^2+||\varepsilon^{s}||^2\right)+
\frac{\tau\kappa^2}{2h^\alpha}\left(||\delta_x^\alpha\varepsilon^s||^2+||\varepsilon^{s+1}||^2\right)\vspace{0.2 cm}\\
&&\displaystyle+\frac{\tau\kappa^2}{6h^\alpha}\lambda_{\max}(A_\alpha)\left(||\varepsilon^{s+1}||^2+||\varepsilon^{s}||^2\right)\vspace{0.2 cm}\\
&\leq&\displaystyle2\left(\frac{\tau\nu^2}{12}+\frac{1}{\tau}\right)\left(||\varepsilon^{s+1}||^2+||\varepsilon^{s}||^2\right)+
\frac{\tau\kappa^2}{2h^\alpha}\left(||\delta_x^\alpha\varepsilon^s||^2+||\varepsilon^{s+1}||^2\right)\vspace{0.2 cm}\\
&&\displaystyle+\frac{2^{2+\alpha}\tau\kappa^2}{6\pi h^\alpha}\left(||\varepsilon^{s+1}||^2+||\varepsilon^{s}||^2\right).
\end{array}
$$
And
$$
\begin{array}{lll}
\displaystyle
W^s&=&\displaystyle\left(\left(\left(\frac{\tau\nu^2}{12}+\frac{1}{\tau}\right)\left({\delta_x^0+\frac{\alpha}{24}\delta_x^2}\right)-\frac{\tau\kappa^2}{6h^\alpha}\delta_x^\alpha\right)
\delta_t\varepsilon^{s+\frac{1}{2}},\delta_t\varepsilon^{s+\frac{1}{2}}\right)+
\frac{\tau\kappa^2}{h^\alpha}\left(\delta_x^\alpha\varepsilon^s,\varepsilon^{s+1}\right)\vspace{0.2 cm}\\
&&\displaystyle+\frac{\tau\kappa^2}{4h^\alpha}\left(\delta_x^\alpha\delta_t\varepsilon^{s+\frac{1}{2}},\delta_t\varepsilon^{s+\frac{1}{2}}\right)\vspace{0.2 cm}\\
&=&\displaystyle\left(\left(\left(\frac{\tau\nu^2}{12}+\frac{1}{\tau}\right)\left({\delta_x^0+\frac{\alpha}{24}\delta_x^2}\right)-\frac{\tau\kappa^2}{6h^\alpha}\delta_x^\alpha\right)
\delta_t\varepsilon^{s+\frac{1}{2}},\delta_t\varepsilon^{s+\frac{1}{2}}\right)+
\frac{\tau\kappa^2}{h^\alpha}\left(\widetilde{\delta}_x^\alpha\varepsilon^s,\widetilde{\delta}_x^\alpha\varepsilon^{s+1}\right)\vspace{0.2 cm}\\
&&\displaystyle+\frac{\tau\kappa^2}{4h^\alpha}\left(\widetilde{\delta}_x^\alpha\delta_t\varepsilon^{s+\frac{1}{2}},\widetilde{\delta}_x^\alpha\delta_t\varepsilon^{s+\frac{1}{2}}\right)\vspace{0.2 cm}\\
&=&\displaystyle\left(\left(\left(\frac{\tau\nu^2}{12}+\frac{1}{\tau}\right)\left({\delta_x^0+\frac{\alpha}{24}\delta_x^2}\right)-\frac{\tau\kappa^2}{6h^\alpha}\delta_x^\alpha\right)
\delta_t\varepsilon^{s+\frac{1}{2}},\delta_t\varepsilon^{s+\frac{1}{2}}\right)+
\frac{\tau\kappa^2}{h^\alpha}\left(\widetilde{\delta}_x^\alpha\varepsilon^{s+\frac{1}{2}},\widetilde{\delta}_x^\alpha\varepsilon^{s+\frac{1}{2}}\right)\vspace{0.2 cm}\\
&&\displaystyle-\frac{\tau\kappa^2}{h^\alpha}\left[\left(\widetilde{\delta}_x^\alpha\varepsilon^{s+\frac{1}{2}},\widetilde{\delta}_x^\alpha\varepsilon^{s+\frac{1}{2}}\right)
-\left(\widetilde{\delta}_x^\alpha\varepsilon^s,\widetilde{\delta}_x^\alpha\varepsilon^{s+1}\right)
\right]
+\frac{\tau\kappa^2}{4h^\alpha}\left(\widetilde{\delta}_x^\alpha\delta_t\varepsilon^{s+\frac{1}{2}},\widetilde{\delta}_x^\alpha\delta_t\varepsilon^{s+\frac{1}{2}}\right)\vspace{0.2 cm}\\
&=&\displaystyle\left(\left(\left(\frac{\tau\nu^2}{12}+\frac{1}{\tau}\right)\left({\delta_x^0+\frac{\alpha}{24}\delta_x^2}\right)-\frac{\tau\kappa^2}{6h^\alpha}\delta_x^\alpha\right)
\delta_t\varepsilon^{s+\frac{1}{2}},\delta_t\varepsilon^{s+\frac{1}{2}}\right)+
\frac{\tau\kappa^2}{h^\alpha}\left({\delta}_x^\alpha\varepsilon^{s+\frac{1}{2}},\varepsilon^{s+\frac{1}{2}}\right)\vspace{0.2 cm}\\
&\geq&\displaystyle\left(\frac{\tau\nu^2}{12}+\frac{1}{\tau}\right)\left(\left(\left({\delta_x^0+\frac{\alpha}{24}\delta_x^2}\right)-\frac{\tau\kappa^2}{6h^\alpha}
\left(\frac{\tau\nu^2}{12}+\frac{1}{\tau}\right)^{-1}\delta_x^\alpha\right)
\delta_t\varepsilon^{s+\frac{1}{2}},\delta_t\varepsilon^{s+\frac{1}{2}}\right)\vspace{0.2 cm}\\&&\displaystyle+
\frac{\tau\kappa^2}{h^\alpha}\lambda_{\min}(A_\alpha)||\varepsilon^{s+\frac{1}{2}}||^2\vspace{0.2 cm}\\
&\geq&\displaystyle\lambda_{\min}(G_\alpha)\left(\frac{\tau\nu^2}{12}+\frac{1}{\tau}\right)||\delta_t\varepsilon^{s+\frac{1}{2}}||^2+
\frac{\tau\kappa^2}{h^\alpha}\lambda_{\min}(A_\alpha)||\varepsilon^{s+\frac{1}{2}}||^2
\vspace{0.2 cm}\\
&\geq&\displaystyle
\lambda_{\min}(G_\alpha)\left(\frac{\tau\nu^2}{12}+\frac{1}{\tau}\right)\left\|\delta_t\varepsilon^{s+\frac{1}{2}}\right\|^2,
\end{array}
$$
Hence, one has
$$
\begin{array}{lll}
\displaystyle
||\delta_t\varepsilon^{s+\frac{1}{2}}||^2&\leq&\displaystyle\frac{1}{\lambda_{\min}(G_\alpha)}\left[
2\left(||\varepsilon^{1}||^2+||\varepsilon^{0}||^2\right)\right.\vspace{0.2 cm}\\
&&\displaystyle\left.+
\frac{6-\alpha}{4g_0^{(\alpha)}}\left(||\varepsilon^{1}||^2+||\delta_x^\alpha\varepsilon^0||^2\right)
+\frac{2^{\alpha}(6-\alpha)}{3\pi g_0^{(\alpha)}}\left(||\varepsilon^{1}||^2+||\varepsilon^{0}||^2\right)
\right],
\end{array}
$$
under the condition (21), again note that $||\varepsilon^0||=||\rho^0||$ and
$
||\varepsilon^1||\leq\frac{6}{6-\alpha}||\rho^1||
$, then
$$
\begin{array}{lll}
\displaystyle
||\delta_t\varepsilon^{s+\frac{1}{2}}||^2&\leq&\displaystyle\frac{1}{\lambda_{\min}(G_\alpha)}\left[
\left(2+\frac{2^{\alpha}(6-\alpha)}{3\pi g_0^{(\alpha)}}\right)||\rho^{0}||^2
+\frac{(6-\alpha)}{4 g_0^{(\alpha)}}||\delta_x^\alpha\rho^0||^2\right.\vspace{0.2 cm}\\&&\displaystyle+
\left.\frac{36}{(6-\alpha)^2}
\left(2+\frac{(6-\alpha)}{4 g_0^{(\alpha)}}+\frac{2^{\alpha}(6-\alpha)}{3\pi g_0^{(\alpha)}}\right)
||\rho^{1}||^2\right]
.
\end{array}
$$
This ends the proof.
\QEDA

\subsection{Convergence analysis}

{\it{\bf Theorem 4}. Under the condition (21), finite
difference scheme (20) is convergent with order $\mathcal{O}(\tau^4+h^4)$.}

{\bf Proof.} Let $e_j^s=u(x_j,t_s)-u_j^s$, $j=1,2,\ldots,M-1, s=0,1,\ldots,N$. From the above analysis, we obtain the following error system:

$$\left\{
\begin{array}{lll} \displaystyle
\displaystyle \left({\delta_x^0+\frac{\alpha}{24}\delta_x^2}\right)J_te_j^s+\frac{\kappa^2}{h^\alpha}
H_t\delta_x^{\alpha}e_{j}^s=\left({\delta_x^0+\frac{\alpha}{24}\delta_x^2}\right)H_tR_j^s,\vspace{0.2 cm}\\ \hspace{2.3 cm}\;\;j=1,\ldots,M-1,\;s=2,\ldots,N-1,\vspace{0.2 cm}\\
  \displaystyle
  e_j^0=0,
\;\;j=1,\ldots,M-1,\vspace{0.2 cm}\\
   \displaystyle
\left({\delta_x^0+\frac{\alpha}{24}\delta_x^2}\right)e_j^1+
\frac{\kappa^2\tau^2}{6h^\alpha}\delta_x^{\alpha}e_{j}^1=R_j^1,\;\;j=1,2,\ldots, M-1.
\vspace{0.2 cm}\\
\displaystyle
e_0^s=e_M^s=0,
\;
\;s=0,1,\ldots,N.
\end{array}\right.\eqno(25)
$$

Firstly, taking the inner product of the third
equation of (25) with $e_j^1$ gets
$$
\begin{array}{lll}
\displaystyle
\left(\left({\delta_x^0+\frac{\alpha}{24}\delta_x^2}\right)e^1,e^1\right)+
\frac{\kappa^2\tau^2}{6h^\alpha}\left(\delta_x^{\alpha}e^1,e^1\right)=\left(R^1,e^1\right).
\end{array}
$$
By using Lemmas 4 and 6 for the above equation, one has
$$
\begin{array}{lll}
\displaystyle
\left(1-\frac{\alpha}{6}\right)||e^1||\leq||R^1||\leq c_1\sqrt{L-l}(\tau^4+h^4),
\end{array}
$$
i.e.,
$$
\begin{array}{lll}
\displaystyle
||e^1||\leq \frac{6c_1}{6-\alpha}\sqrt{L-l}(\tau^4+h^4).
\end{array}
$$

Next, taking the inner product of the first
equation of (25) with $\mu_t\delta_t e_j^s$, one gets
$$
\begin{array}{lll}
\displaystyle
\left(\left({\delta_x^0+\frac{\alpha}{24}\delta_x^2}\right)J_te^s,\mu_t\delta_te^s\right)
+\frac{\kappa^2}{h^\alpha}\left(
H_t\delta_x^{\alpha}e^s,\mu_t\delta_te^s\right)\vspace{0.2 cm}\\\displaystyle=\left(\left({\delta_x^0+\frac{\alpha}{24}\delta_x^2}\right)H_tR^s,\mu_t\delta_te^s\right).
\end{array}\eqno(26)
$$

Denote
$$
\begin{array}{lll}
\displaystyle
E^s&=&\displaystyle\left(\frac{\tau\nu^2}{12}+\frac{1}{\tau}\right)\left(\left({\delta_x^0+\frac{\alpha}{24}\delta_x^2}\right)\delta_te^{s+\frac{1}{2}},\delta_te^{s+\frac{1}{2}}\right)+
\frac{\tau\kappa^2}{h^\alpha}\left(\delta_x^\alpha e^s,e^{s+1}\right)\vspace{0.2 cm}\\
&&\displaystyle+\frac{\tau\kappa^2}{12h^\alpha}\left(\delta_x^\alpha\delta_t e^{s+\frac{1}{2}},\delta_t e^{s+\frac{1}{2}}\right).
\end{array}
$$
Similar to the stability analysis, equ. (26) can be rewritten as
$$
\begin{array}{lll}
\displaystyle
E^s+\frac{\nu}{2 }\left(\left({\delta_x^0+\frac{\alpha}{24}\delta_x^2}\right)\mu_t\delta_t e^s,\mu_t\delta_t e^s\right)
+\frac{\nu\tau^2\kappa^2}{24h^\alpha}\left(\delta_x^\alpha\mu_t\delta_t e^s,\mu_t\delta_t e^s\right)\vspace{0.2 cm}\\ \displaystyle =E^{s-1}+\tau\left(\left({\delta_x^0+\frac{\alpha}{24}\delta_x^2}\right)H_tR^s,\mu_t\delta_te^s\right).
\end{array}\eqno(27)
$$

For the right-hand term of equ. (27), we have the following estimate by using Cauchy-Schwarz inequality and Lemma 9,
$$
\begin{array}{lll}
&&\displaystyle
\tau\left(\left({\delta_x^0+\frac{\alpha}{24}\delta_x^2}\right)H_tR^s,\mu_t\delta_te^s\right)\vspace{0.2 cm}\\&\leq&\displaystyle\tau
\left(||R^s||+\frac{1}{12}||\delta_t^2R^s||+\frac{\nu\tau}{2}||\mu_t\delta_tR^s||\right)
\left(\left\|\delta_te^{s+\frac{1}{2}}\right\|+\left\|\delta_te^{s-\frac{1}{2}}\right\|\right)

\vspace{0.2 cm}\\&\leq&\displaystyle
\frac{\tau}{\lambda_{\min}(G_\alpha)}\left(||R^s||^2+\frac{1}{12}||\delta_t^2R^s||^2+\frac{\nu\tau}{2}||\mu_t\delta_tR^s||^2\right)\vspace{0.2 cm}\\&&\displaystyle
+\tau\lambda_{\min}(G_\alpha)
\left(\frac{13}{24}+\frac{\nu\tau}{4}\right)\left(\left\|\delta_te^{s+\frac{1}{2}}\right\|^2+\left\|\delta_te^{s-\frac{1}{2}}\right\|^2\right)

\vspace{0.2 cm}\\&\leq&\displaystyle
\frac{\tau}{\lambda_{\min}(G_\alpha)}\left(2||R^s||^2+\left(\frac{1}{4}+\nu\tau\right)\left(||R^{s+1}||^2+||R^{s-1}||^2\right)\right)\vspace{0.2 cm}\\&&\displaystyle
+\tau\lambda_{\min}(G_\alpha)
\left(\frac{13}{24}+\frac{\nu\tau}{4}\right)\left(\left\|\delta_te^{s+\frac{1}{2}}\right\|^2+\left\|\delta_te^{s-\frac{1}{2}}\right\|^2\right).
\end{array}
$$
Now equ. (27) can be rewritten as
$$
\begin{array}{lll}
\displaystyle
E^s&\leq&\displaystyle E^{s-1}+\tau\lambda_{\min}(G_\alpha)
\left(\frac{13}{24}+\frac{\nu\tau}{4}\right)\left(\left\|\delta_te^{s+\frac{1}{2}}\right\|^2+\left\|\delta_te^{s-\frac{1}{2}}\right\|^2\right)
\vspace{0.2 cm}\\
&&\displaystyle+\frac{\tau}{\lambda_{\min}(G_\alpha)}\left(2||R^s||^2+\left(\frac{1}{4}+\nu\tau\right)\left(||R^{s+1}||^2+||R^{s-1}||^2\right)\right).
\end{array}\eqno(28)
$$
Similarly, one also has
$$
\begin{array}{lll}
\displaystyle
E^s\geq\lambda_{\min}(G_\alpha)\left(\frac{\tau\nu^2}{12}+\frac{1}{\tau}\right)\left\|\delta_te^{s+\frac{1}{2}}\right\|^2
\geq\frac{\sqrt{3}}{3}\lambda_{\min}(G_\alpha)\left\|\delta_te^{s+\frac{1}{2}}\right\|^2.
\end{array}\eqno(29)
$$
Combining (28) with (29) yields,
$$
\begin{array}{lll}
\displaystyle
E^s&\leq&\displaystyle E^{s-1}+\sqrt{3}\tau\left(\frac{13}{24}+\frac{\nu\tau}{4}\right)\left(E^s+E^{s-1}\right)\vspace{0.2 cm}\\
&&\displaystyle
+\frac{\tau}{\lambda_{\min}(G_\alpha)}\left(2||R^s||^2+\left(\frac{1}{4}+\nu\tau\right)\left(||R^{s+1}||^2+||R^{s-1}||^2\right)\right).
\end{array}
$$
If $\tau\left(\frac{13}{24}+\frac{\nu\tau}{4}\right)\leq\frac{\sqrt{3}}{9}$, we have
$$
\begin{array}{lll}
\displaystyle
E^s&\leq&\displaystyle\left(1+6\sqrt{3}\tau\left(\frac{13}{24}+\frac{\nu\tau}{4}\right)\right) E^{s-1}\vspace{0.2 cm}\\&&\displaystyle
+\frac{3\tau}{2\lambda_{\min}(G_\alpha)}\left(2||R^s||^2+\left(\frac{1}{4}+\nu\tau\right)\left(||R^{s+1}||^2+||R^{s-1}||^2\right)\right).
\end{array}
$$
Utilizing the Gronwall inequality, one has
$$
\begin{array}{lll}
\displaystyle
E^s&\leq&\displaystyle e^{6\sqrt{3}\left(\frac{13}{24}+\frac{\nu\tau}{4}\right)s\tau} \left[E^{0}
+\frac{3\tau}{2\lambda_{\min}(G_\alpha)}\sum_{k=1}^{s}\left(2||R^k||^2\right.\right.\vspace{0.2 cm}\\&&\displaystyle
\left.\left.+\left(\frac{1}{4}+\nu\tau\right)\left(||R^{k+1}||^2+||R^{k-1}||^2\right)\right)\right]\vspace{0.2 cm}\\
&\leq&\displaystyle e^{6\sqrt{3}\left(\frac{13}{24}+\frac{\nu}{4}\right)T} \left[E^{0}
+\frac{3T(L-l)}{2\lambda_{\min}(G_\alpha)}\left(2+\left(\frac{1}{4}+\nu\right)\left(c_1^2+2c_2^2\right)\right)\left(\tau^4+h^4\right)^2\right].
\end{array}
$$
The above inequality can be
changed into,
$$
\begin{array}{lll}
\displaystyle
||\delta_te^{s+\frac{1}{2}}||^2&\leq&\displaystyle\frac{(L-l)e^{6\sqrt{3}\left(\frac{13}{24}+\frac{\nu}{4}\right)T}}{\lambda_{\min}(G_\alpha)}\left[
\frac{36c_1^2(L-l)}{(6-\alpha)^2}
\left(2+\frac{(6-\alpha)}{4 g_0^{(\alpha)}}+\frac{2^{\alpha}(6-\alpha)}{3\pi g_0^{(\alpha)}}\right)\right.\vspace{0.2 cm}\\&&\displaystyle
\left.
+\frac{3T}{2\lambda_{\min}(G_\alpha)}\left(2+\left(\frac{1}{4}+\nu\right)\left(c_1^2+2c_2^2\right)\right)\right]\left(\tau^4+h^4\right)^2
,
\end{array}
$$
i.e.,
$$
\begin{array}{lll}
\displaystyle
||\delta_te^{s+\frac{1}{2}}||\leq C\left(\tau^4+h^4\right),
\end{array}
$$
where $C=C_1(C_2+C_3)$,
$$
\begin{array}{lll}
\displaystyle
C_1=\exp\left({3\sqrt{3}\left(\frac{13}{24}+\frac{\nu}{4}\right)T}\right)\sqrt{\frac{(L-l)}{\lambda_{\min}(G_\alpha)}}\;,
\end{array}
$$
$$
\begin{array}{lll}
\displaystyle
C_2=\frac{6c_1}{6-\alpha}\sqrt{\left(L-l\right)
\left(2+\frac{(6-\alpha)}{4 g_0^{(\alpha)}}+\frac{2^{\alpha}(6-\alpha)}{3\pi g_0^{(\alpha)}}\right)},
\end{array}
$$
and
$$
\begin{array}{lll}
\displaystyle
C_3=\sqrt{\frac{3T}{2\lambda_{\min}(G_\alpha)}\left(2+\left(\frac{1}{4}+\nu\right)\left(c_1^2+2c_2^2\right)\right)}\;.
\end{array}
$$
The proof is thus completed.
\QEDA

 \section{Numerical examples}
 \quad
 In this section we present some numerical experiments.

\textbf{Example 1.} Consider
 function
$u_n(x)=x^{2n}(1-x)^{2n}$, $x\in[0,1]$, $n=2,3,4,5.$

To illustrate the convergence orders of the numerical schemes (11)--(14), we
compute the absolute error $E(h)$ by
$$
\begin{array}{lll}
\displaystyle
E(h)=\left|u\left(x_{\frac{1}{2h}}\right)-u_{\frac{1}{2h}}\right|,
\end{array}
$$
where $u\left(x_{\frac{1}{2h}}\right)$ represents the exact solution and
$u_{\frac{1}{2h}}$ is the numerical solution of function $u(x)$ at point $x_{\frac{1}{2h}}$.
In other words, we compute the absolute error and convergence orders at point $x=\frac{1}{2}$
by different stepsize $h$ due to $x_{\frac{1}{2h}}=\left(\frac{1}{2h}\right)h=\frac{1}{2}$.

We calculate the experimental convergence order (denoted by ECO) by
$$
\begin{array}{lll}
\displaystyle
\textrm{ECO}=\left. \log\left(\frac{E(h_1)}{E(h_2)}\right) \middle/ \log\left(\frac{h_1}{h_2}\right) \right..
\end{array}
$$

 From Definition 1, we know that the Riesz derivative of the above function is
$$
\begin{array}{lll}
\displaystyle \frac{\partial^\alpha
u(x)}{\partial{|x|^\alpha}}=-\frac{1}{2\cos\left(\pi\alpha/2\right)}
\sum
\limits_{\ell=0}^{{2n}}(-1)^{\ell}\frac{(2n)!({2n}+l)!\left[x^{{2n}+\ell-\alpha}+(1-x)^{{2n}+\ell-\alpha}\right]}{\ell!({2n}-l)!\mathrm{\Gamma}({2n}+\ell+1-\alpha)}
.
\end{array}
$$

We numerically solve $u_n(x)$ by using numerical schemes (11)--(14) with various values of
 $h$ and $\alpha$. The $E(h)$
and ECO are shown in Tables \ref{tab.1}-\ref{tab.4}.
 From these tables, it is clear that the experimental
 orders are in line with the theoretical orders.

\begin{table}[!htbp]
\begin{center}
\caption{ The absolute error $E(h)$ and the experimental convergence order (ECO) of function $u_2(x)$ by
numerical scheme (11).}\label{tab.1}
 \vspace{0.3 cm}
 \begin{footnotesize}
\begin{tabular}{c c c c c c }\hline \\
  $\alpha$ & $h$&  $E(h)$&ECO
 \vspace{0.2 cm} \\ \hline
   $1.1 $& $\frac{1}{20}$ &1.985528e-006 & ---\vspace{0.1 cm}\\
  $$  & $\frac{1}{40}$&  1.247417e-007  &  3.9925\vspace{0.1 cm}\\
   $$  & $\frac{1}{80}$&7.806460e-009 &  3.9981\vspace{0.1 cm}\\
  $$& $\frac{1}{160}$ &4.880592e-010  &  3.9995\vspace{0.1 cm}\\
   $$& $\frac{1}{320}$ &  3.050456e-011  &  4.0000\vspace{0.1 cm}\\\hline
   $1.3 $& $\frac{1}{20}$ &3.418165e-006 &      ---\vspace{0.1 cm}\\
  $$  & $\frac{1}{40}$& 2.147477e-007   &    3.9925\vspace{0.1 cm}\\
   $$  & $\frac{1}{80}$& 1.343913e-008  &   3.9981\vspace{0.1 cm}\\
  $$& $\frac{1}{160}$ &8.402173e-010   & 3.9995\vspace{0.1 cm}\\
   $$& $\frac{1}{320}$ & 5.251979e-011  &  	3.9998\vspace{0.1 cm}\\ \hline
  $1.5$& $\frac{1}{20}$ &5.712995e-006 & ---\vspace{0.1 cm}\\
  $$  & $\frac{1}{40}$&  	3.588944e-007	  & 3.9926\vspace{0.1 cm}\\
   $$  & $\frac{1}{80}$&2.245955e-008	  &     3.9982\vspace{0.1 cm}\\
  $$& $\frac{1}{160}$ &1.404168e-009    & 3.9995\vspace{0.1 cm}\\
  $$& $\frac{1}{320}$ & 8.777425e-011	  &  3.9998\vspace{0.1 cm}\\\hline
  $1.7 $& $\frac{1}{20}$ &9.316621e-006  &      ---\vspace{0.1 cm}\\
  $$  & $\frac{1}{40}$&  5.851888e-007    &  3.9928\vspace{0.1 cm}\\
   $$  & $\frac{1}{80}$&3.661963e-008   &   3.9982\vspace{0.1 cm}\\
  $$& $\frac{1}{160}$ & 2.289431e-009	     & 3.9996\vspace{0.1 cm}\\
   $$& $\frac{1}{320}$ & 1.431370e-010	 &   3.9995\vspace{0.1 cm}\\ \hline
  $1.9 $& $\frac{1}{20}$ &1.486627e-005	  &      ---\vspace{0.1 cm}\\
  $$  & $\frac{1}{40}$& 9.335587e-007 &        3.9932\vspace{0.1 cm}\\
   $$  & $\frac{1}{80}$& 5.841643e-008      &  3.9983 \vspace{0.1 cm}\\
  $$& $\frac{1}{160}$ &3.652104e-009	   &      3.9996\vspace{0.1 cm}\\
   $$& $\frac{1}{320}$ &  2.284402e-010  &       	3.9988\vspace{0.1 cm}\\
 \hline
\end{tabular}
 \end{footnotesize}
 \end{center}
 \end{table}

\begin{table}[!htbp]
\begin{center}
\caption{ The absolute error $E(h)$ and the experimental convergence order (ECO) of function $u_3(x)$ by
numerical scheme (12).}\label{tab.2}
 \vspace{0.3 cm}
 \begin{footnotesize}
\begin{tabular}{c c c c c c }\hline \\
  $\alpha$ & $h$&   $E(h)$&ECO
 \vspace{0.2 cm} \\ \hline
   $1.1 $& $\frac{1}{20}$ &3.120201e-008 & ---\vspace{0.1 cm}\\
  $$  & $\frac{1}{24}$&  1.057512e-008 &  5.9345\vspace{0.1 cm}\\
   $$  & $\frac{1}{28}$&4.223802e-009	 & 5.9537\vspace{0.1 cm}\\
  $$& $\frac{1}{32}$ &1.904343e-009	 & 	5.9656\vspace{0.1 cm}\\
   $$& $\frac{1}{36}$ & 9.422903e-010  &  	5.9735\vspace{0.1 cm}\\\hline
   $1.3 $& $\frac{1}{20}$ &6.008620e-008	 &      ---\vspace{0.1 cm}\\
  $$  & $\frac{1}{24}$& 2.036744e-008  &   5.9337\vspace{0.1 cm}\\
   $$  & $\frac{1}{28}$& 8.135715e-009 &  5.9531\vspace{0.1 cm}\\
  $$& $\frac{1}{32}$ &3.668326e-009	  & 5.9651\vspace{0.1 cm}\\
   $$& $\frac{1}{36}$ &1.815230e-009	  &  		5.9730\vspace{0.1 cm}\\ \hline
  $1.5$& $\frac{1}{20}$ &1.123916e-007 & ---\vspace{0.1 cm}\\
  $$  & $\frac{1}{24}$&  	3.810009e-008	  & 5.9333\vspace{0.1 cm}\\
   $$  & $\frac{1}{28}$&1.521978e-008	  &    5.9528\vspace{0.1 cm}\\
  $$& $\frac{1}{32}$ &	6.862760e-009	  &5.9648\vspace{0.1 cm}\\
  $$& $\frac{1}{36}$ &3.396074e-009		  & 5.9728\vspace{0.1 cm}\\\hline
  $1.7 $& $\frac{1}{20}$ &2.053203e-007 &      ---\vspace{0.1 cm}\\
  $$  & $\frac{1}{24}$&  6.960225e-008    & 5.9333\vspace{0.1 cm}\\
   $$  & $\frac{1}{28}$&2.780411e-008	&   	5.9527\vspace{0.1 cm}\\
  $$& $\frac{1}{32}$ & 1.253731e-008		     & 5.9647\vspace{0.1 cm}\\
   $$& $\frac{1}{36}$ & 	6.204232e-009		 &  5.9727\vspace{0.1 cm}\\ \hline
  $1.9 $& $\frac{1}{20}$ &3.675466e-007	  &      ---\vspace{0.1 cm}\\
  $$  & $\frac{1}{24}$& 1.245861e-007 &        5.9338\vspace{0.1 cm}\\
   $$  & $\frac{1}{28}$& 4.976658e-008	    &  5.9530\vspace{0.1 cm}\\
  $$& $\frac{1}{32}$ &2.244006e-008		   &     5.9649\vspace{0.1 cm}\\
   $$& $\frac{1}{36}$ &  1.110461e-008  &       	5.9727\vspace{0.1 cm}\\
 \hline
\end{tabular}
 \end{footnotesize}
 \end{center}
 \end{table}

\begin{table}[!htbp]
\begin{center}
\caption{ The absolute error $E(h)$ and the experimental convergence order (ECO) of function $u_4(x)$ by
numerical scheme (13).}\label{tab.3}
 \vspace{0.3 cm}
 \begin{footnotesize}
\begin{tabular}{c c c c c c }\hline \\
  $\alpha$ & $h$&  $E(h)$&ECO
 \vspace{0.2 cm} \\ \hline
   $1.1 $& $\frac{1}{30}$ &3.442344e-011 & ---\vspace{0.1 cm}\\
  $$  & $\frac{1}{34}$&  1.279869e-011 &  7.9048\vspace{0.1 cm}\\
   $$  & $\frac{1}{38}$&5.303531e-012	 &7.9206\vspace{0.1 cm}\\
  $$& $\frac{1}{42}$ &2.397829e-012 & 	7.9315\vspace{0.1 cm}\\
   $$& $\frac{1}{46}$ & 1.164521e-012	  &  	7.9393\vspace{0.1 cm}\\\hline
   $1.3 $& $\frac{1}{30}$ &	7.195825e-011	 &      ---\vspace{0.1 cm}\\
  $$  & $\frac{1}{34}$& 2.675551e-011  &  7.9044\vspace{0.1 cm}\\
   $$  & $\frac{1}{38}$& 1.108608e-011 & 7.9213\vspace{0.1 cm}\\
  $$& $\frac{1}{42}$ &5.011612e-012  & 7.9327\vspace{0.1 cm}\\
   $$& $\frac{1}{46}$ &2.433728e-012	  &  	7.9402\vspace{0.1 cm}\\ \hline
  $1.5$& $\frac{1}{30}$ &1.459781e-010 & ---\vspace{0.1 cm}\\
  $$  & $\frac{1}{34}$&  	5.427806e-011	  &7.9044\vspace{0.1 cm}\\
   $$  & $\frac{1}{38}$&2.248828e-011	  &   7.9220\vspace{0.1 cm}\\
  $$& $\frac{1}{42}$ &	1.016534e-011		  &7.9335\vspace{0.1 cm}\\
  $$& $\frac{1}{46}$ &4.936491e-012		  & 7.9401\vspace{0.1 cm}\\\hline
  $1.7 $& $\frac{1}{30}$ &2.889640e-010 &      ---\vspace{0.1 cm}\\
  $$  & $\frac{1}{34}$&  1.074284e-010  &7.9055\vspace{0.1 cm}\\
   $$  & $\frac{1}{38}$&4.449502e-011	&   	7.9249\vspace{0.1 cm}\\
  $$& $\frac{1}{42}$ & 2.010185e-011		     & 7.9390\vspace{0.1 cm}\\
   $$& $\frac{1}{46}$ & 	9.752953e-012		 & 7.9502\vspace{0.1 cm}\\ \hline
  $1.9 $& $\frac{1}{30}$ &5.601931e-010	  &      ---\vspace{0.1 cm}\\
  $$  & $\frac{1}{34}$& 	2.082444e-010 &       7.9062\vspace{0.1 cm}\\
   $$  & $\frac{1}{38}$& 8.624075e-011	    &  7.9260\vspace{0.1 cm}\\
  $$& $\frac{1}{42}$ &3.895758e-011	   &     7.9401\vspace{0.1 cm}\\
   $$& $\frac{1}{46}$ & 	1.890060e-011	  &      7.9506\vspace{0.1 cm}\\
 \hline
\end{tabular}
 \end{footnotesize}
 \end{center}
 \end{table}

\begin{table}[!htbp]
\begin{center}
\caption{ The absolute error $E(h)$ and the experimental convergence order (ECO) of function $u_5(x)$ by
numerical scheme (14).}\label{tab.4}
 \vspace{0.3 cm}
 \begin{footnotesize}
\begin{tabular}{c c c c c c }\hline \\
  $\alpha$ & $h$&  $E(h)$&ECO
 \vspace{0.2 cm} \\ \hline
   $1.1 $& $\frac{1}{30}$ &2.669378e-012 & ---\vspace{0.1 cm}\\
  $$  & $\frac{1}{34}$&  6.799466e-013 &  9.5568\vspace{0.1 cm}\\
   $$  & $\frac{1}{38}$&2.057061e-013		 &9.5521\vspace{0.1 cm}\\
  $$& $\frac{1}{42}$ &7.247668e-014 & 	9.3790\vspace{0.1 cm}\\
   $$& $\frac{1}{46}$ & 2.983513e-014	  &  	8.8684\vspace{0.1 cm}\\\hline
   $1.3 $& $\frac{1}{30}$ &6.076893e-012	 &      ---\vspace{0.1 cm}\\
  $$  & $\frac{1}{34}$&1.530707e-012  &  9.6349\vspace{0.1 cm}\\
   $$  & $\frac{1}{38}$&4.563785e-013	& 9.6687\vspace{0.1 cm}\\
  $$& $\frac{1}{42}$ &1.566398e-013& 9.6145\vspace{0.1 cm}\\
   $$& $\frac{1}{46}$ &6.127250e-014	  &  	9.3784\vspace{0.1 cm}\\ \hline
  $1.5$& $\frac{1}{30}$ &1.316180e-011& ---\vspace{0.1 cm}\\
  $$  & $\frac{1}{34}$&  	3.294738e-012	  &9.6784\vspace{0.1 cm}\\
   $$  & $\frac{1}{38}$&	9.747428e-013	  &   9.7306\vspace{0.1 cm}\\
  $$& $\frac{1}{42}$ &	3.301910e-013		  &9.7325\vspace{0.1 cm}\\
  $$& $\frac{1}{46}$ &1.257671e-013		  & 	9.6443\vspace{0.1 cm}\\\hline
  $1.7 $& $\frac{1}{30}$ &2.756366e-011 &      ---\vspace{0.1 cm}\\
  $$  & $\frac{1}{34}$& 6.859984e-012  &9.7190\vspace{0.1 cm}\\
   $$  & $\frac{1}{38}$&2.007171e-012		&   		9.8190\vspace{0.1 cm}\\
  $$& $\frac{1}{42}$ & 6.624227e-013			     & 9.9669\vspace{0.1 cm}\\
   $$& $\frac{1}{46}$ & 2.368039e-013		 &10.2781\vspace{0.1 cm}\\ \hline
  $1.9 $& $\frac{1}{30}$ &5.640743e-011	 &      ---\vspace{0.1 cm}\\
  $$  & $\frac{1}{34}$& 	1.401455e-011 &      9.7309\vspace{0.1 cm}\\
   $$  & $\frac{1}{38}$&4.097488e-012	    &  	9.8250\vspace{0.1 cm}\\
  $$& $\frac{1}{42}$ &1.354007e-012   &    9.9555\vspace{0.1 cm}\\
   $$& $\frac{1}{46}$ & 		4.868649e-013		  &      	10.2198\vspace{0.1 cm}\\
 \hline
\end{tabular}
 \end{footnotesize}
 \end{center}
 \end{table}

{\bf Example 2.} Consider function $u(x)=x(1-x)$ for $x\in[0,1]$, which doesn't meet the conditions of the Theorem 2 and
Remark 1. We numerically compute this function by using 4th-order formula (11), the convergence order
is list in Table \ref{tab.5}. From the table, one can see that the expected convergence order 4 in (11) cann't be obtained,
which implies that the suitable smoothness assumptions in Theorem 2 or Remark 1 are needed.

\begin{table}[!htbp]
\begin{center}
\caption{ The numerical results of the Example 2 by using 4th-order fractional-compact
formula (11).}\label{tab.5}
 \vspace{0.3 cm}
 \begin{footnotesize}
\begin{tabular}{c c c c c c }\hline \\
  $\alpha$ & $h$&  $E(h)$&ECO
 \vspace{0.2 cm} \\ \hline
   $1.1 $& $\frac{1}{10}$ &5.900848e-04	 & ---\vspace{0.1 cm}\\
  $$  & $\frac{1}{20}$&  1.470732e-04	  &  2.0044\vspace{0.1 cm}\\
   $$  & $\frac{1}{40}$&3.674047e-05 &  2.0011\vspace{0.1 cm}\\
  $$& $\frac{1}{80}$ &9.183380e-06	  &  2.0003\vspace{0.1 cm}\\
   $$& $\frac{1}{160}$ &  2.295736e-06	  & 2.0001\vspace{0.1 cm}\\\hline
   $1.3 $& $\frac{1}{10}$ &6.828118e-04 &      ---\vspace{0.1 cm}\\
  $$  & $\frac{1}{20}$& 1.699820e-04	   &   2.0061\vspace{0.1 cm}\\
   $$  & $\frac{1}{40}$& 4.245071e-05	  &   2.0015\vspace{0.1 cm}\\
  $$& $\frac{1}{80}$ &1.060988e-05	   & 2.0004\vspace{0.1 cm}\\
   $$& $\frac{1}{160}$ & 2.652296e-06	  &  	2.0001\vspace{0.1 cm}\\ \hline
  $1.5$& $\frac{1}{10}$ &7.105214e-04	 & ---\vspace{0.1 cm}\\
  $$  & $\frac{1}{20}$&  	1.766378e-04	&	2.0081\vspace{0.1 cm}\\
   $$  & $\frac{1}{40}$&  4.409782e-05		&2.0020\vspace{0.1 cm}\\
  $$& $\frac{1}{80}$ &1.102061e-05		&2.0005\vspace{0.1 cm}\\
  $$& $\frac{1}{160}$ & 2.754912e-06	&	2.0001\vspace{0.1 cm}\\\hline
  $1.7 $& $\frac{1}{10}$ &6.101600e-04	  &      ---\vspace{0.1 cm}\\
  $$  & $\frac{1}{20}$&  1.514533e-04	&	2.0103\vspace{0.1 cm}\\
   $$  & $\frac{1}{40}$&3.779590e-05	&	2.0026\vspace{0.1 cm}\\
  $$& $\frac{1}{80}$ & 9.444768e-06		&2.0006\vspace{0.1 cm}\\
   $$& $\frac{1}{160}$ & 2.360928e-06	&	2.0002\vspace{0.1 cm}\\ \hline
  $1.9 $& $\frac{1}{10}$ &2.863584e-04		  &      ---\vspace{0.1 cm}\\
  $$  & $\frac{1}{20}$& 7.095794e-05	&	2.0128\vspace{0.1 cm}\\
   $$  & $\frac{1}{40}$& 1.770032e-05	&	2.0032 \vspace{0.1 cm}\\
  $$& $\frac{1}{80}$ &4.422636e-06		&  2.0008\vspace{0.1 cm}\\
   $$& $\frac{1}{160}$ & 1.105502e-06	&	2.0002\vspace{0.1 cm}\\
 \hline
\end{tabular}
 \end{footnotesize}
 \end{center}
 \end{table}

{\bf Example 3.} Consider the following equation
$$
\begin{array}{rrr} \displaystyle
 \frac{\partial{{}^2u(x,t)}}{\partial{t^2}}+\frac{\partial{{}u(x,t)}}{\partial{t}}
 =\frac{\partial^\alpha u(x,t)}{\partial{|x|^\alpha}}
  +f(x,t),\;\;1<\alpha<2,\\\,\,\,0<x<1,\,\,\,0< t\leq 1$$,
\end{array}
$$
where
$$
\begin{array}{lll}
\displaystyle f(x,t)&=&\displaystyle 2(2t^2+t+1)\exp(t^2)
x^6(1-x)^6+\frac{1}{2}\exp(t^2)\sec\left(\frac{\pi}{2}\alpha\right)\times
\\&&\displaystyle\left\{\frac{\mathrm{\Gamma}(7)}{\mathrm{\Gamma}(7-\alpha)}
\left[x^{6-\alpha}+(1-x)^{6-\alpha}\right]-\frac{6\mathrm{\Gamma}(8)}{\mathrm{\Gamma}(8-\alpha)}
\left[x^{7-\alpha}+(1-x)^{7-\alpha}\right]\right.\vspace{0.2
cm}\\&&\displaystyle
+\frac{15\mathrm{\Gamma}(9)}{\mathrm{\Gamma}(9-\alpha)}
\left[x^{8-\alpha}+(1-x)^{8-\alpha}\right]
\displaystyle-\frac{20\mathrm{\Gamma}(10)}{\mathrm{\Gamma}(10-\alpha)}
\left[x^{9-\alpha}+(1-x)^{9-\alpha}\right]\vspace{0.2 cm}\\&&\displaystyle
+\frac{15\mathrm{\Gamma}(11)}{\mathrm{\Gamma}(11-\alpha)}
\left[x^{10-\alpha}+(1-x)^{10-\alpha}\right] -\frac{6\mathrm{\Gamma}(12)}{\mathrm{\Gamma}(12-\alpha)}
\left[x^{11-\alpha}+(1-x)^{11-\alpha}\right] \vspace{0.2 cm}\\&&
\displaystyle\left.
+\frac{\mathrm{\Gamma}(13)}{\mathrm{\Gamma}(13-\alpha)}
\left[x^{12-\alpha}+(1-x)^{12-\alpha}\right] \right\}.
\end{array}
$$

Its exact solution is given by $\displaystyle u(x, t) =x^6(1-x)^6\exp(t)$, which satisfy the necessary initial
and boundary value conditions.

Here, the absolute error is calculated by the formula
$$
\begin{array}{lll}
\displaystyle
E(\tau,h)=\max_{0\leq j\leq M,0\leq s\leq N}\left|u\left(x_j,t_s\right)-u_j^s\right|,
\end{array}
$$
the spatial experimental convergence order (denoted by SECO) and temporal experimental
 convergence order (denoted by TECO) are both computed by the following formula
$$
\begin{array}{lll}
\displaystyle
\textrm{SECO}=\textrm{TECO}=\left. \log_2\left(\frac{E(2\tau,2h)}{E(\tau,h)}\right) \right.,
\end{array}
$$
respectively.

\begin{figure}
\centering
 \includegraphics[width=11 cm]{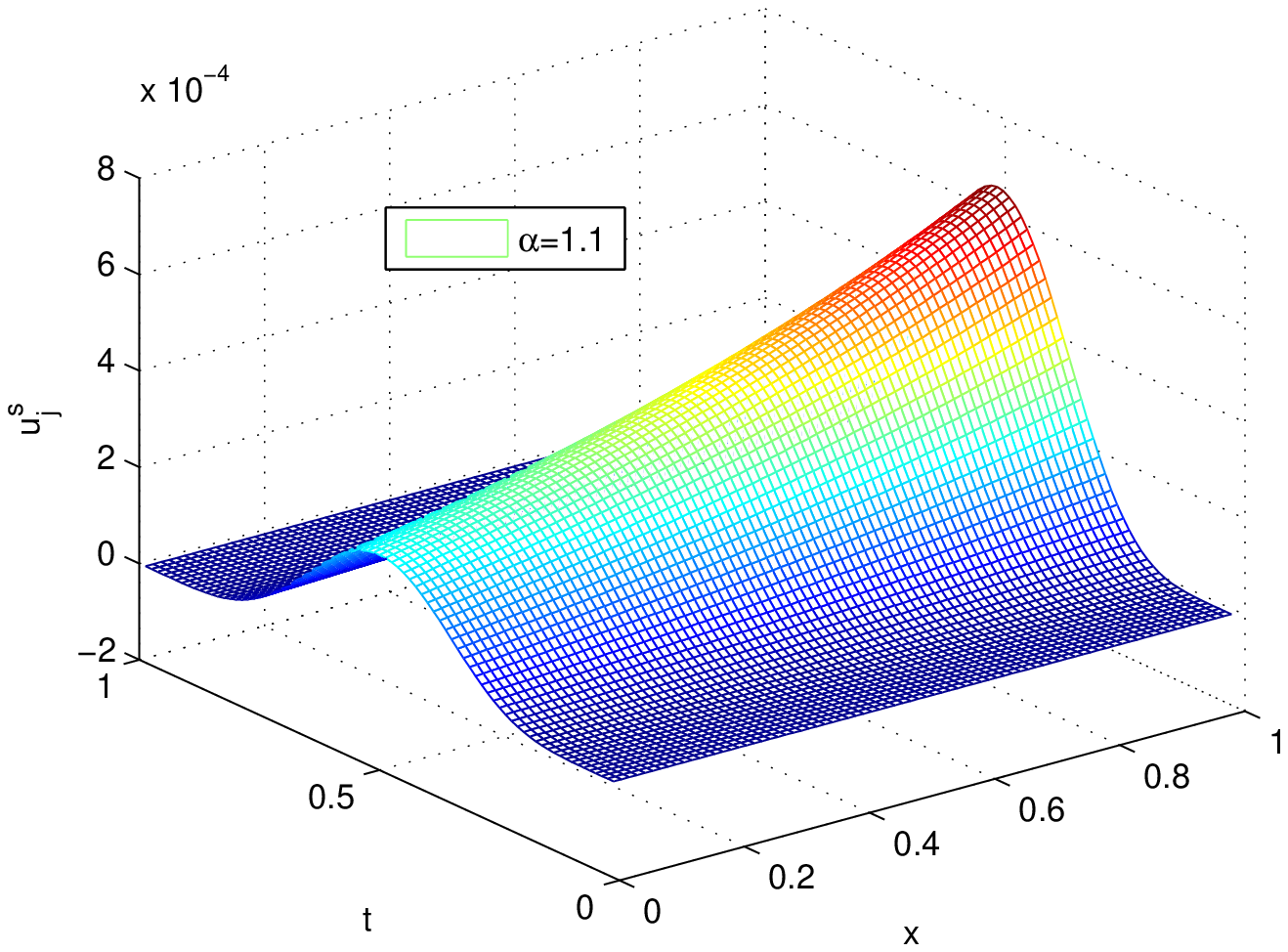}\\
  \caption{The numerical solution surface by difference scheme (20) with $\alpha=1.1$ and $\tau=h=\frac{1}{80}$.}
  \label{fig.4}
\end{figure}

The numerical solution surface by the proposed finite difference scheme (20)
with different $h$, $\tau$ and $\alpha$ (which satisfy the stability
condition (21))
are given in Figures \ref{fig.4}-\ref{fig.7}. The absolute error, temporal and spatial
 experimental convergence orders are listed in Table \ref{tab.6}.
  It can be concluded from the table and the figures that the convergence
order of the finite difference scheme (20) is $\mathcal{O}(\tau^4+h^4)$, which are in
line with the theoretical analysis.
\begin{table}\renewcommand\arraystretch{1.2}
 \begin{center}
 \caption{ The absolute error $E(\tau,h)$, temporal and spatial experimental convergence
orders (TECO and SECO) of the Example 3 by difference scheme (20).}\label{tab.6}
 \vspace{0.3 cm}
 \begin{footnotesize}
\begin{tabular}{c c c c c  }\hline
   & \;\;\;\;  &\;\;\;\;\;\;&
   \\
  $\alpha$ &\;\;&  \;\;$E(\tau,h)$\;\;\;\; & TECO\;(SECO)
  \\
  \hline \vspace{0.1 cm}
  $1.2 $& $h=\frac{1}{4},\tau=\frac{1}{4}$ &    6.629148e-005&  ---
\\ \vspace{0.1 cm}
  $$  & $h=\frac{1}{8},\tau=\frac{1}{8}$&         3.374742e-006 &   4.2960\\ \vspace{0.1 cm}
  $$& $h=\frac{1}{16},\tau=\frac{1}{16}$ &       1.987545e-007&    4.0857\\ \vspace{0.1 cm}
$$  & $h=\frac{1}{32},\tau=\frac{1}{32}$&        1.216316e-008 &     4.0304\\\vspace{0.1 cm}
  $$& $h=\frac{1}{64},\tau=\frac{1}{64}$ &        7.535224e-010&  4.0127 \\ \vspace{0.1 cm}
$$  & $h=\frac{1}{128},\tau=\frac{1}{128}$&        4.689844e-011 &   4.0060\\
 \hline\vspace{0.1 cm}
   $1.4 $& $h=\frac{1}{4},\tau=\frac{1}{4}$ &  7.814928e-005&  ---
\\ \vspace{0.1 cm}
 $$  & $h=\frac{1}{8},\tau=\frac{1}{8}$&        3.293314e-006&    4.5686\\\vspace{0.1 cm}
  $$& $h=\frac{1}{16},\tau=\frac{1}{16}$ &      1.906135e-007&   4.1108 \\\vspace{0.1 cm}
$$  & $h=\frac{1}{32},\tau=\frac{1}{32}$&        1.160778e-008&   4.0375\\\vspace{0.1 cm}
  $$& $h=\frac{1}{64},\tau=\frac{1}{64}$ &      7.178332e-010&   4.0153 \\\vspace{0.1 cm}
$$  & $h=\frac{1}{128},\tau=\frac{1}{128}$&      4.465339e-011 &  4.0068\\
 \hline\vspace{0.1 cm}
  $1.6 $& $h=\frac{1}{4},\tau=\frac{1}{4}$&      9.026191e-005
&  ---\;\;&    ---
\\ \vspace{0.1 cm}
 $$  & $h=\frac{1}{8},\tau=\frac{1}{8}$&        3.701763e-006 &    4.6078\\\vspace{0.1 cm}
  $$& $h=\frac{1}{16},\tau=\frac{1}{16}$ &       2.350544e-007&  3.9771 \\\vspace{0.1 cm}
$$  & $h=\frac{1}{32},\tau=\frac{1}{32}$&        1.451435e-008&   4.0174\\\vspace{0.1 cm}
  $$& $h=\frac{1}{64},\tau=\frac{1}{64}$ &        9.043217e-010&   4.0045\\\vspace{0.1 cm}
$$  & $h=\frac{1}{128},\tau=\frac{1}{128}$&        5.651652e-011 &   4.0001\\
 \hline\vspace{0.1 cm}
   $1.8 $& $h=\frac{1}{4},\tau=\frac{1}{4}$ &     1.027098e-004&  ---
\\ \vspace{0.1 cm}
  $$  & $h=\frac{1}{8},\tau=\frac{1}{8}$&       4.552933e-006 &  4.4956\\\vspace{0.1 cm}
  $$& $h=\frac{1}{16},\tau=\frac{1}{16}$ &       2.329253e-007&    4.2889 \\\vspace{0.1 cm}
$$  & $h=\frac{1}{32},\tau=\frac{1}{32}$&        1.358572e-008 & 4.0997\\\vspace{0.1 cm}
  $$& $h=\frac{1}{64},\tau=\frac{1}{64}$ &      8.403641e-010&   4.0149
 \\\vspace{0.1 cm}
$$  & $h=\frac{1}{128},\tau=\frac{1}{128}$&        5.231801e-011 &  4.0056\\
 \hline
\end{tabular}
 \end{footnotesize}
 \end{center}
 \end{table}

\begin{figure}
\centering
 \includegraphics[width=11 cm]{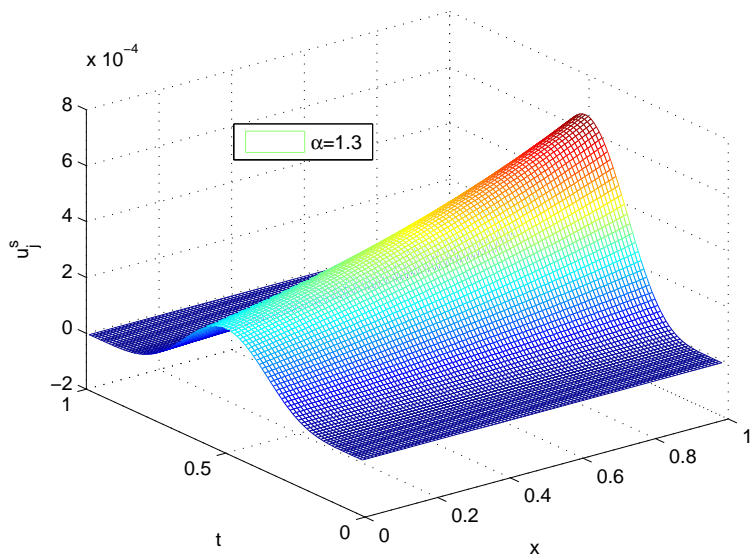}\\
  \caption{The numerical solution surface by difference scheme (20) with $\alpha=1.3$ and $\tau=h=\frac{1}{100}$.}
  \label{fig.5}
\end{figure}

\section{Conclusions}
\quad\;
 In this work, we firstly derive even order
 fractional-compact numerical differential formulas for Riesz derivatives. Then a 4th-order
 numerical formula has been applied to solving the Riesz spatial fractional
 telegraph equation. Numerical results show that
 the experimental convergence orders are close to the theoretical convergence orders and confirm
 the efficiency of the new method developed in the current paper.
 These fractional-compact numerical differential formulas and the techniques can be also applied to $m$-dimensional $(m\geq2)$
  Riesz spatial differential equations.\\

\textbf{{Appendix I}}

Firstly, we list two lemmas as follows:

{\it{\bf Lemma I}. (\cite{KJC}) The following inequalities hold:
 $$\begin{array}{lll}
 \mathrm{(i)}\;\;1-x<\exp\left(-x\right)\;\;\textmd{for} \;\; 0<x<1.\vspace{0.1 cm}\\
\displaystyle \mathrm{(ii)}\;\;1-x>\exp\left(-2x\right)\;\;\textmd{for} \;\; 0<x\leq0.7968.
\end{array}
$$}

{\it{\bf Lemma II}. If functions $f(x)$ and $g(x)$ satisfy the conditions: $f(x)<0$ is an increasing
function and $g(x)>0$ is a decreasing
function, then $f(x)g(x)$ is an increasing function.}

The proof is easy so is omitted here.

Now, we give the detailed proof for Theorem 1. Note that ${g}_{k}^{(\alpha)}={g}_{-k}^{(\alpha)}$, so we only consider the case $k\geq0$.

{\bf Proof.} (\romannumeral1)  From the recursive relation
$$\begin{array}{rrr}
\displaystyle
{g}_{_k}^{(\alpha)}=\left(1-\frac{\alpha+1}{\frac{\alpha}{2}+k}\right){g}_{_{k-1}}^{(\alpha)},
\end{array}\eqno(I1)
$$
we easily obtain
$$\begin{array}{lll}
\displaystyle
\left|{g}_{_k}^{(\alpha)}\right|&<&\displaystyle\exp\left(-\frac{\alpha+1}{\frac{\alpha}{2}+k}\right)\left|{g}_{_{k-1}}^{(\alpha)}\right|<
\exp\left(-\frac{\alpha+1}{\frac{\alpha}{2}+k}\right)\exp\left(-\frac{\alpha+1}{\frac{\alpha}{2}+k-1}\right)
\left|{g}_{_{k-2}}^{(\alpha)}\right|\vspace{0.2 cm}
\\
&<&\displaystyle\cdots<\exp\left(-\frac{\alpha+1}{\frac{\alpha}{2}+k}\right)\exp\left(-\frac{\alpha+1}{\frac{\alpha}{2}+k-1}\right)\ldots
\exp\left(-\frac{\alpha+1}{\frac{\alpha}{2}+3}\right)
\left|{g}_{_{2}}^{(\alpha)}\right|
\vspace{0.2 cm}\\&=&\displaystyle
\exp\left(-(\alpha+1)\sum_{n=3}^{k}\frac{1}{\frac{\alpha}{2}+n}\right)\left|{g}_{_{2}}^{(\alpha)}\right|
\end{array}
$$
by using Lemma I.

Because function $\frac{1}{\frac{\alpha}{2}+x}$ is a decreasing function for $x>0$, then
$$\begin{array}{rrr}
\displaystyle
\sum_{n=3}^{k}\frac{1}{\frac{\alpha}{2}+n}>\int_{3}^{k+1} \frac{1}{\frac{\alpha}{2}+x}dx=\ln\frac{\frac{\alpha}{2}+k+1}{\frac{\alpha}{2}+3},
\end{array}
$$
so,
$$\begin{array}{rrr}
\displaystyle
\left|{g}_{_k}^{(\alpha)}\right|<
\exp\left(-(\alpha+1)\ln\frac{\frac{\alpha}{2}+k+1}{\frac{\alpha}{2}+3}\right)\left|{g}_{_{2}}^{(\alpha)}\right|
=-S(\alpha)\left(\frac{{\alpha}+6}{{\alpha}+2(k+1)}\right)^{\alpha+1},\;\;k\geq3,
\end{array}
$$
where $S(\alpha)=\frac{\mathrm{\Gamma}(\alpha+1)}{\mathrm{\Gamma}\left(\frac{\alpha}{2}-1\right)
\mathrm{\Gamma}\left(\frac{\alpha}{2}+3\right)}$.

Similarly, we also have the following estimation by equ. $(I1)$ and Lemma I,
$$\begin{array}{rrr}
\displaystyle
\left|{g}_{_k}^{(\alpha)}\right|>\exp\left(-\frac{2(\alpha+1)}{\frac{\alpha}{2}+k}\right)\left|{g}_{_{k-1}}^{(\alpha)}\right|>
\exp\left(-2(\alpha+1)\sum_{n=3}^{k}\frac{1}{\frac{\alpha}{2}+n}\right)\left|{g}_{_{2}}^{(\alpha)}\right|.
\end{array}
$$
Notice
$$\begin{array}{rrr}
\displaystyle
\sum_{n=3}^{k}\frac{1}{\frac{\alpha}{2}+n}<\int_{2}^{k} \frac{1}{\frac{\alpha}{2}+x}dx=\ln\frac{\frac{\alpha}{2}+k}{\frac{\alpha}{2}+2},
\end{array}
$$
one has
$$\begin{array}{rrr}
\displaystyle
\left|{g}_{_k}^{(\alpha)}\right|>
\exp\left(-2(\alpha+1)\ln\frac{\frac{\alpha}{2}+k}{\frac{\alpha}{2}+2}\right)\left|{g}_{_{2}}^{(\alpha)}\right|
=-S(\alpha)\left(\frac{{\alpha}+4}{{\alpha}+2k}\right)^{2(\alpha+1)},\;\;k\geq3.
\end{array}
$$

(\romannumeral2) Since
$$\begin{array}{lll}
\displaystyle
\sum_{k=n}^{m}\frac{1}{(\alpha+2k)^{2(\alpha+1)}}&>&\displaystyle\int_{n}^{m+1} \frac{1}{(\alpha+2x)^{2(\alpha+1)}}dx\\
&=&\displaystyle\frac{1}{2(2\alpha+1)}\left[\left(\alpha+2n\right)^{-2\alpha-1}-\left(\alpha+2m+2\right)^{-2\alpha-1}\right],
\end{array}
$$
and
$$\begin{array}{lll}
\displaystyle
\sum_{k=n}^{m}\frac{1}{(\alpha+2k+2)^{\alpha+1}}&<&\displaystyle\int_{n-1}^{m} \frac{1}{(\alpha+2x+2)^{\alpha+1}}dx\\
&=&\displaystyle\frac{1}{2\alpha}\left[\left(\alpha+2n\right)^{-\alpha}-\left(\alpha+2m+2\right)^{-\alpha}\right],
\end{array}
$$
one gets,
$$\begin{array}{lll}
\displaystyle
\sum_{k=n}^{m}\left|g_{_k}^{(\alpha)}\right|&>&\displaystyle-S(\alpha)\left(\alpha+4\right)^{2(\alpha+1)}
\sum_{k=n}^{m}\frac{1}{(\alpha+2k)^{2(\alpha+1)}}\vspace{0.2 cm}\\
&>&\displaystyle
\displaystyle\frac{\left(\alpha+4\right)^{2(\alpha+1)}}{2(2\alpha+1)}P_1(m,n,\alpha),
\end{array}\eqno(I2)
$$
where $P_1(m,n,\alpha)=-S(\alpha)
\left[\left(\alpha+2n\right)^{-2\alpha-1}-\left(\alpha+2m+2\right)^{-2\alpha-1}\right]$, $n\geq3$,
and
$$\begin{array}{lll}
\displaystyle
\sum_{k=n}^{m}\left|g_k^{(\alpha)}\right|&<&\displaystyle-S(\alpha)\left(\alpha+6\right)^{\alpha+1}
\sum_{k=n}^{m}\frac{1}{(\alpha+2(k+1))^{\alpha+1}}\\\vspace{0.2 cm}&<&\displaystyle
\displaystyle\frac{\left(\alpha+6\right)^{\alpha+1}}{2\alpha}P_2(m,n,\alpha),
\end{array}\eqno(I3)
$$
where $P_2(m,n,\alpha)=-S(\alpha)
\left[\left(\alpha+2n\right)^{-\alpha}-\left(\alpha+2m+2\right)^{-\alpha}\right]$, $n\geq3$.

(\romannumeral3) Just letting $m=\infty$ in $(I2)$ and $(I3)$ gives the proof.

(\romannumeral4) By using the following formula \cite{EMOT}
$$\begin{array}{lll}
\displaystyle
\int_{0}^{\frac{\pi}{2}}\left(\cos x\right)^{\alpha}\cos\left(\beta x\right)\mathrm{d}x
=\frac{\pi}{2^{1+\alpha}}\frac{\mathrm{\Gamma}(1+\alpha)}{\mathrm{\Gamma}\left
(1+\frac{\alpha+\beta}{2}\right)\mathrm{\Gamma}\left(1+\frac{\alpha-\beta}{2}\right)},\;\;\;\mathrm{Re}(\alpha)>-1,
\end{array}
$$
one can gets
$$\begin{array}{lll}
\displaystyle
g_{_0}^{(\alpha)}=\frac{\mathrm{\Gamma}(1+\alpha)}{\mathrm{\Gamma}\left
(1+\frac{\alpha}{2}\right)\mathrm{\Gamma}\left(1+\frac{\alpha}{2}\right)}=
\frac{2^{1+\alpha}}{\pi}\int_{0}^{\frac{\pi}{2}}\left(\cos x\right)^{\alpha}\mathrm{d}x
\leq\frac{2^{1+\alpha}}{\pi}\int_{0}^{\frac{\pi}{2}}\cos x\mathrm{d}x=\frac{2^{1+\alpha}}{\pi}
\end{array}
$$
for $1<\alpha<2$.

In addition, $g_{_0}^{(\alpha)}$ is usually represented by using the beta function $B(\cdot,\cdot)$, i.e.,
$$\begin{array}{lll}
\displaystyle
g_{_0}^{(\alpha)}=\frac{\mathrm{\Gamma}(1+\alpha)}{\mathrm{\Gamma}\left
(1+\frac{\alpha}{2}\right)\mathrm{\Gamma}\left(1+\frac{\alpha}{2}\right)}=
\frac{\mathrm{\Gamma}(2+\alpha)}{(1+\alpha)\mathrm{\Gamma}\left
(1+\frac{\alpha}{2}\right)\mathrm{\Gamma}\left(1+\frac{\alpha}{2}\right)}
=\frac{1}{(1+\alpha) B\left(\frac{2+\alpha}{2},\frac{2+\alpha}{2}\right)},
\end{array}
$$
and in which
$$\begin{array}{lll}
\displaystyle
 B\left(\frac{\alpha}{2},\frac{\alpha}{2}\right)
 =2\int_{0}^{\frac{\pi}{2}}\left(\sin x\right)^{\alpha-1}\left(\cos x\right)^{\alpha-1}\mathrm{d}x=
 \frac{1}{2^{\alpha-2}}\int_{0}^{\frac{\pi}{2}}\left(\sin 2x\right)^{\alpha-1}\mathrm{d}x
 \leq\frac{\pi}{2^{\alpha-1}},
\end{array}
$$
so, we have
$$\begin{array}{lll}
\displaystyle
g_{_0}^{(\alpha)}\geq\frac{2^{1+\alpha}}{\pi(1+\alpha)}.
\end{array}
$$

From \cite{SLAT} we known that
$$\begin{array}{lll}
\displaystyle
\sum\limits_{k=-\infty,k\neq0}^
{\infty}\left|\right.{g}_{_k}^{(\alpha)}\left|\right.={g}_{_0}^{(\alpha)},
\end{array}
$$
that is,
$$\begin{array}{lll}
\displaystyle
\frac{2^{1+\alpha}}{(1+\alpha)\pi}\leq\sum\limits_{k=-\infty,k\neq0}^
{\infty}\left|\right.{g}_{_k}^{(\alpha)}\left|\right.\leq\frac{2^{1+\alpha}}{\pi}.
\end{array}
$$

(\romannumeral5) When $k=0$, then $g_{_0}^{(\alpha)}=\frac{\mathrm{\Gamma}(\alpha+1)}{\mathrm{\Gamma}^2\left(\frac{\alpha}{2}+1\right)}$,
and $\ln g_{_0}^{(\alpha)}=\ln\mathrm{\Gamma}(\alpha+1)-2\ln\mathrm{\Gamma}\left(\frac{\alpha}{2}+1\right)$.
$$\begin{array}{lll}
\displaystyle
\frac{ g_0^{(\alpha)^{'}}}{g_0^{(\alpha)}}&=&\displaystyle\psi(\alpha+1)-\psi\left(\frac{\alpha}{2}+1\right)\vspace{0.2 cm}\\
&=&\displaystyle\int_{0}^{\infty}\frac{\exp(-t)-\exp(-(\alpha+1)t)}{1-\exp(-t)}\mathrm{d}t
-\int_{0}^{\infty}\frac{\exp(-t)-\exp\left(-\left(\frac{\alpha}{2}+1\right)t\right)}{1-\exp\left(-t\right)}\mathrm{d}t
\vspace{0.2 cm}\\
&=&\displaystyle\int_{0}^{\infty}\frac{\exp(-(\alpha+1)t)\left[
\exp\left(\frac{\alpha}{2}t\right)-1
\right]}{1-\exp(-t)}\mathrm{d}t
\vspace{0.2 cm}\\
&=&\displaystyle
\int_{0}^{\infty}\frac{\exp(-(\alpha+1)t)
}{1-\exp(-t)}\left(\sum_{n=0}^{\infty}\frac{\left(\frac{\alpha}{2}\right)^n}{n!}t^n\right)\mathrm{d}t>0,
\end{array}
$$
where $g_0^{(\alpha)^{'}}=\frac{\mathrm{d}g_0^{(\alpha)}}{\mathrm{d} \alpha}$. Since $g_0^{(\alpha)}>0$, $g_0^{(\alpha)^{'}}>0$, that is, $g_0^{(\alpha)}$
is an increasing function with respect to $\alpha\in(1,2)$.

 When $k=1$, $g_1^{(\alpha)}=-\frac{\mathrm{\Gamma}(\alpha+1)}
 {\mathrm{\Gamma}\left(\frac{\alpha}{2}\right)\mathrm{\Gamma}\left(\frac{\alpha}{2}+2\right)}$, similarly to case $k=0$, we have

$$\begin{array}{lll}
\displaystyle
\frac{\mathrm{d}}{\mathrm{d}\alpha}\ln\left(-g_1^{(\alpha)}\right)
&=&\displaystyle\frac{ g_1^{(\alpha)^{'}}}{g_1^{(\alpha)}}=\psi(\alpha+1)-\frac{1}{2}\psi\left(\frac{\alpha}{2}\right)
-\frac{1}{2}\psi\left(\frac{\alpha}{2}+2\right)
\vspace{0.2 cm}\\
&=&\displaystyle
\int_{0}^{\infty}\frac{\frac{1}{2}\exp\left(\frac{-\alpha}{2}t\right)+\frac{1}{2}\exp\left(-\left(\frac{\alpha}{2}+2\right)t\right)-
\exp(-(\alpha+1)t)
}{1-\exp(-t)}\mathrm{d}t
\vspace{0.2 cm}\\
&\geq&\displaystyle
\int_{0}^{\infty}\frac{\exp\left(-\left(\frac{\alpha}{2}+1\right)t\right)-
\exp(-(\alpha+1)t)
}{1-\exp(-t)}\mathrm{d}t>0.
\end{array}
$$
Because of $g_1^{(\alpha)}<0$, then $g_1^{(\alpha)^{'}}<0$, i.e, $g_1^{(\alpha)}$
is a decreasing function.

When $k=2$,
$$\begin{array}{lll}
\displaystyle
g_2^{(\alpha)}=\frac{\mathrm{\Gamma}(\alpha+1)}{\mathrm{\Gamma}\left(\frac{\alpha}{2}+3\right)\mathrm{\Gamma}\left(\frac{\alpha}{2}-1\right)}
=-\frac{2\left(1-\frac{\alpha}{2}\right)\mathrm{\Gamma}(\alpha)}{\left(1+\frac{\alpha}{2}\right)\left(2+\frac{\alpha}{2}
\right)\mathrm{\Gamma}^2\left(\frac{\alpha}{2}\right)}.
\end{array}
$$
According to another definition of $\mathrm{\Gamma}(\alpha)$ \cite{SKM}:
$$\begin{array}{lll}
\displaystyle
\mathrm{\Gamma}(\alpha)=\frac{1}{\alpha}\prod\limits_{n=1}^{\infty}\frac{\left(1+\frac{1}{n}\right)^\alpha}{1+\frac{\alpha}{n}},
\end{array}
$$
it is easily known that
$$\begin{array}{lll}
\displaystyle
\frac{\mathrm{\Gamma}(\alpha)}{\mathrm{\Gamma}^2\left(\frac{\alpha}{2}\right)}
=
\frac{\alpha}{4}\prod\limits_{n=1}^{\infty}\frac{\left(1+\frac{\alpha}{2n}\right)^2}{1+\frac{\alpha}{n}}.
\end{array}
$$

Substituting the above equation into $g_2^{(\alpha)}
$ yields
$$\begin{array}{lll}
\displaystyle
g_2^{(\alpha)}
=-\frac{\alpha(2-\alpha)}{(2+\alpha)(4+\alpha)}\prod\limits_{n=1}^{\infty}
\frac{\left(1+\frac{\alpha}{2n}\right)^2}{1+\frac{\alpha}{n}}.
\end{array}
$$
$$\begin{array}{lll}
\displaystyle
\frac{g_2^{(\alpha)^{'}}}{g_2^{(\alpha)}}=
\frac{\mathrm{d}}{\mathrm{d}\alpha}\ln\left(-g_2^{(\alpha)}\right)=\frac{1}{2}\sum_{n=1}^{\infty}\frac{\alpha}
{\left(n+\frac{\alpha}{2}\right)(n+\alpha)}
-\left(\frac{1}{2-\alpha}-\frac{1}{\alpha}+\frac{1}{2+\alpha}+\frac{1}{4+\alpha}\right).
\end{array}
$$

Let
$$\begin{array}{lll}
\displaystyle
S(\alpha)=
\frac{1}{2}\sum_{n=1}^{\infty}\frac{1}{\left(n+\frac{\alpha}{2}\right)(n+\alpha)},
\end{array}
$$
and
$$\begin{array}{lll}
\displaystyle
f_1(\alpha)=\frac{1}{\alpha}
\left(\frac{1}{2-\alpha}-\frac{1}{\alpha}+\frac{1}{2+\alpha}+\frac{1}{4+\alpha}\right).
\end{array}
$$
Then
$$\begin{array}{lll}
\displaystyle
S(\alpha)<
\sum_{n=1}^{\infty}\frac{2}{(2n+1)^2}=\sum_{n=1}^{\infty}\frac{2}{(2n-1)^2}-2=\frac{\pi^2}{4}-2,
\end{array}
$$
and
$$\begin{array}{lll}
\displaystyle
f_1^{'}(\alpha)=\frac{8\left[3\alpha^5+10\alpha^3+12+8(2-\alpha)+4(2\alpha-3)^2\right]}{\alpha^3(2-\alpha)^2(2+\alpha)^2(4+\alpha)^2}
>0,
\end{array}
$$
that is to say $f_1(\alpha)$ is an increasing function and its minimum value is $f_1(\alpha)_{\min}=f_1(1)=\frac{8}{15}.$

So,
$$\begin{array}{lll}
\displaystyle
\frac{g_2^{(\alpha)^{'}}}{g_2^{(\alpha)}}=
\alpha(S(\alpha)-f_1(\alpha))<\alpha\left(\frac{\pi^2}{4}-2
-\frac{8}{15}\right)<\alpha\left(\frac{3.15^2}{4}-2
-\frac{8}{15}\right)=-\frac{253}{4800}\alpha.
\end{array}
$$
Since $g_2^{(\alpha)}<0$, it follows that $g_2^{(\alpha)^{'}}>0$, i.e.,
$g_2^{(\alpha)}$ is an increasing function with respect to $\alpha\in(1,2)$.

When $k=3$,
$$\begin{array}{lll}
\displaystyle
g_3^{(\alpha)}
=-\frac{2\left(1-\frac{\alpha}{2}\right)\left(2-\frac{\alpha}{2}\right)\mathrm{\Gamma}(\alpha)}
{\left(1+\frac{\alpha}{2}\right)\left(2+\frac{\alpha}{2}\right)\left(3+\frac{\alpha}{2}\right)\mathrm{\Gamma}^2\left(\frac{\alpha}{2}\right)}
=-\frac{\alpha(2-\alpha)(4-\alpha)}{(2+\alpha)(4+\alpha)(6+\alpha)}
\prod\limits_{n=1}^{\infty}\frac{\left(1+\frac{\alpha}{2n}\right)^2}{1+\frac{\alpha}{n}}
.
\end{array}
$$
$$\begin{array}{lll}
\displaystyle
\frac{g_3^{(\alpha)^{'}}}{g_3^{(\alpha)}}=
\alpha(S(\alpha)-f_2(\alpha)),
\end{array}
$$
where $\displaystyle f_2(\alpha)=\frac{1}{\alpha}\left(\frac{1}{2-\alpha}+\frac{1}{4-\alpha}
-\frac{1}{\alpha}+\frac{1}{2+\alpha}+\frac{1}{4+\alpha}+\frac{1}{6+\alpha}\right)>0.$

So,
$$\begin{array}{lll}
\displaystyle
\frac{g_3^{(\alpha)^{'}}}{g_3^{(\alpha)}}=
\alpha(S(\alpha)-f_2(\alpha))<\alpha\left(\frac{\pi^2}{4}-2
-\frac{106}{105}\right)<\alpha\left(\frac{3.15^2}{4}-2
-\frac{106}{105}\right)=-\frac{485}{917}\alpha.
\end{array}
$$
It immediately follows that $g_3^{(\alpha)}$ is an increasing function.

Finally, we consider the cases $k\geq4$. From \cite{O} one has
$$\begin{array}{lll}
\displaystyle
g_k^{(\alpha)}=\frac{(-1)^k\mathrm{\Gamma}(\alpha+1)}{\mathrm{\Gamma}\left(\frac{\alpha}{2}-k+1\right)\mathrm{\Gamma}\left(\frac{\alpha}{2}+k+1\right)}
=-\frac{\sin\left(\frac{\pi}{2}\alpha\right)}{\pi}\frac{\mathrm{\Gamma}(\alpha+1)\mathrm{\Gamma}
\left(k-\frac{\alpha}{2}\right)}{\mathrm{\Gamma}\left(\frac{\alpha}{2}+k+1\right)}.
\end{array}
$$

Letting
$$\begin{array}{lll}
\displaystyle
z(\alpha)
=-\frac{\sin\left(\frac{\pi}{2}\alpha\right)}{\pi},\;\;
 z_{_k}(\alpha)=
\frac{\mathrm{\Gamma}(\alpha+1)\mathrm{\Gamma}\left(k-\frac{\alpha}{2}\right)}{\mathrm{\Gamma}\left(\frac{\alpha}{2}+k+1\right)},
\end{array}
$$
obviously, $z(\alpha)<0$ is an increasing function with respect to $\alpha$ for $1<\alpha<2$. Meanwhile, we have
$$\begin{array}{lll}
\displaystyle
\ln z_{_k}(\alpha)
=\ln\mathrm{\Gamma}(\alpha+1)+\ln\mathrm{\Gamma}\left(k-\frac{\alpha}{2}\right)-\ln\mathrm{\Gamma}\left(\frac{\alpha}{2}+k+1\right).
\end{array}
$$
Furthermore, it follows that
$$\begin{array}{lll}
\displaystyle
\frac{  z_{_k}(\alpha)^{'}}{  z_{_k}(\alpha)}
&=&\displaystyle
\psi(\alpha+1)-\frac{1}{2}\psi\left(k-\frac{\alpha}{2}\right)-\frac{1}{2}\psi\left(\frac{\alpha}{2}+k+1\right)\vspace{0.2 cm}\\
&=&\displaystyle
\int_{0}^{\infty}\frac{\frac{1}{2}\exp\left(-\left(k-\frac{\alpha}{2}\right)t\right)+\frac{1}{2}\exp\left(-\left(\frac{\alpha}{2}+k+1\right)t\right)-
\exp(-(\alpha+1)t)
}{1-\exp(-t)}\mathrm{d}t
\vspace{0.2 cm}\\
&=&\displaystyle
\int_{0}^{\infty}\frac{\exp\left(-\left(\frac{\alpha}{2}+k+1\right)t\right)\left[\frac{1}{2}
\exp((\alpha+1)t)+\frac{1}{2}-\exp((k-\frac{\alpha}{2})t)\right]
}{1-\exp(-t)}\mathrm{d}t
\vspace{0.2 cm}\\
&=&\displaystyle
\int_{0}^{\infty}\frac{\exp\left(-\left(\frac{\alpha}{2}+k+1\right)t\right)}{1-\exp(-t)}
\left(\sum_{n=1}^{\infty}\frac{\frac{1}{2}(\alpha+1)^n-\left(k-\frac{\alpha}{2}\right)^n}{n!}
t^n\right)\mathrm{d}t.
\end{array}
$$
If $k\geq4$, then $\frac{1}{2}(\alpha+1)^n-\left(k-\frac{\alpha}{2}\right)^n<0$ for $1<\alpha<2$. At this moment, we have
$  z_{_k}(\alpha)^{'}<0$, i.e., $ z_{_k}(\alpha)$ is a decreasing function. So $g_k^{(\alpha)}=z(\alpha) z_{_k}(\alpha)$ is an increasing function
for $k\geq4$ by using Lemma II.

Therefore, all this finishes the proof of Theorem 1.\QEDA
\vspace{0.2 cm}\\

\end{document}